\documentclass[reqno]{amsart}
\usepackage{amssymb,epsf}

\newtheorem{theorem}{Theorem}

\newtheorem{corollary}[theorem]{Corollary}

\theoremstyle{definition}
\newtheorem*{definition}{Definition}
\newtheorem*{conjecture}{Conjecture}

\theoremstyle{remark}

\numberwithin{equation}{section}

\begin{document}

\title[Theta Hypergeometric Integrals]
{Theta hypergeometric integrals}

\author{V.P. Spiridonov}
 \address{Bogoliubov Laboratory of Theoretical Physics, Joint
Institute for Nuclear Research, Dubna, Moscow Region 141980, Russia.}

\thanks{Published in:
{\em Algebra i Analiz (St. Petersburg Math. J.)}
{\bf 15} (2003), No. 6, 161--215. }

\begin{abstract}
A general class of (multiple) hypergeometric type integrals
associated with the Jacobi theta functions is defined. These integrals
are related to theta hypergeometric series via the residue calculus.
In the one variable case, theta function extensions of the Meijer
function are obtained. A number of multiple generalizations
of the elliptic beta integral \cite{spi:elliptic}
associated with the root systems $A_n$ and $C_n$ is described.
Some of the $C_n$-examples were proposed earlier by van
Diejen and the author, but other integrals are new.
An example of the biorthogonality relations associated with the
elliptic beta integrals is considered in detail.
\end{abstract}

\maketitle

\centerline{\em Dedicated to Mizan Rahman}
\medskip

\tableofcontents

\section{Introduction}

Exact integration formulas and
integral representations of functions are important from
various points of view. Such representations serve sometimes as
definitions of functions, but more often they are needed for the
better understanding of properties of functions defined beforehand.
Due to numerous applications (see \cite{aar:special}), the Euler beta integral
\begin{equation}\label{euler}
\int_0^1 x^{\alpha-1}(1-x)^{\beta-1}dx=\frac{\Gamma(\alpha)\Gamma(\beta)}
{\Gamma(\alpha+\beta)}, \qquad \mbox{Re }\alpha,\, \mbox{Re }\beta>0,
\end{equation}
where $\Gamma(z)$ is the standard gamma function,
plays a fundamental role in classical analysis.
Various $q$-generalizations of (\ref{euler}) involving $q$-gamma functions
have been proposed within the theory of basic hypergeometric
series \cite{gas-rah:basic}. Recently, a ``third floor" of the
hierarchy of beta-type integrals, which is related to the elliptic
gamma function, has been discovered in \cite{spi:elliptic}
(in the one variable case) and in \cite{die-spi:elliptic,die-spi:selberg}
(multiple extensions). For a brief review of results in this direction,
see \cite{die-spi:review}. In this paper, we discuss a general class of
hypergeometric type integrals associated with the Jacobi theta
functions and propose several new multiple elliptic beta
integrals admitting exact evaluations.

An infinite hierarchy of multiple gamma functions was proposed by
Barnes long time ago \cite{bar:multiple}:
\begin{equation}
\Gamma_r^{-1}(u;\mathbf{\omega}) =
e^{\sum_{i=0}^r\gamma_{ri}\frac{u^i}{i!}}\;u \prod_{n_1,\ldots,
n_r=0}^\infty{}' \left(1+\frac{u}{\Omega}\right)
e^{\sum_{i=1}^r(-1)^i\frac{u^i}{i\Omega^i}},
\label{m-gamma}\end{equation}
where $\gamma_{ri}$ are some constants analogous
to the Euler constant and $\Omega=n_1\omega_1+\ldots+n_r\omega_r$
(if some of the ratios $\omega_i/\omega_k$ are real, then they must be
positive). The prime in the product sign means that the point
$n_1=\ldots=n_r=0$ is skipped. The function (\ref{m-gamma}) satisfies
a collection of $r$ first order difference equations
$$
\frac{\Gamma_r(u+\omega_j;\mathbf{\omega})}
{\Gamma_r(u;\mathbf{\omega})}=
\frac{1}{\Gamma_{r-1}(u;\mathbf{\omega}(j))},
\qquad j=1,\ldots,r,
$$
where $\mathbf{\omega}(j)=(\omega_1,\ldots,\omega_{j-1},\omega_{j+1},\ldots,
\omega_r)$ and $\Gamma_1(u;\omega)=
\rho(\omega) \omega^{u/\omega}\Gamma(u/\omega)$
for some constant $\rho(\omega)$ (for a brief account of this
function, see also appendix A in \cite{jim-miw:quantum}).

Following Barnes' analysis, in \cite{jac:basic} Jackson has considered
the generalized gamma functions in a slightly different way and proposed
the $q$-gamma function and the elliptic gamma function. We recall the
definition of the latter. Taking two complex variables $q$ and $p$ such that
$|q|, |p|<1$, we compose the following (convergent) Jackson double infinite
product:
\begin{equation}
(z;q,p)_\infty=\prod_{j,k=0}^\infty(1-zq^jp^k).
\label{d-product}\end{equation}
Two first order $q$- and $p$-difference equations for this product
\begin{equation}
\frac{(z;q,p)_\infty}{(qz;q,p)_\infty}=(z;p)_\infty, \qquad
\frac{(z;q,p)_\infty}{(pz;q,p)_\infty}=(z;q)_\infty,
\label{eq-1}\end{equation}
where $(z;p)_\infty=\prod_{k=0}^\infty(1-zp^k)$,
are of major importance. Replacing $z$ by $pz^{-1}$ in the first
equation and by $qz^{-1}$ in the second, we get
\begin{equation}
\frac{(qp(qz)^{-1};q,p)_\infty}{(qpz^{-1};q,p)_\infty}
=(pz^{-1};p)_\infty,\qquad
\frac{(qp(pz)^{-1};q,p)_\infty}{(qpz^{-1};q,p)_\infty}
=(qz^{-1};q)_\infty.
\label{eq-2}\end{equation}

We define a theta function as follows:
\begin{equation}
\theta(z;p)=(z;p)_\infty(pz^{-1};p)_\infty.
\label{theta}\end{equation}
It is related to the standard Jacobi $\theta_1$-function
\cite{whi-wat:course} in a simple way
\begin{eqnarray} \nonumber
&& \theta_1(u;\sigma,\tau)=
-i\sum_{n=-\infty}^\infty (-1)^np^{(2n+1)^2/8}q^{(n+1/2)u}
\\ \nonumber && \makebox[2em]{}
=p^{1/8} iq^{-u/2}\: (p;p)_\infty\: \theta(q^{u};p), \quad
u\in\mathbb{C},
\label{theta1}\end{eqnarray}
where we assume that $q=e^{2\pi i\sigma}, p=e^{2\pi i\tau}$.
Sometimes, for brevity, it is convenient to drop $q$ and $p$ or
the modular parameters $\sigma$ and $\tau$ in the notations
for theta functions and elliptic gamma functions, as well as for the
elliptic analogs of shifted factorials to be defined below.

The function $\theta_1(u)$ is entire, odd,
$\theta_1(-u)=-\theta_1(u)$, and doubly quasiperiodic
\begin{equation}
\theta_1(u+\sigma^{-1}) = -\theta_1(u), \qquad
\theta_1(u+\tau\sigma^{-1}) = -e^{-\pi i\tau-2\pi i\sigma u}
\theta_1(u).
\label{quasi}\end{equation}
These transformation properties of the $\theta_1$-function
are extensively used in our formalism. For the
$\theta(z;p)$ function, they take the form
\begin{equation}\label{fun-rel}
\theta(pz;p)=\theta(z^{-1};p)=-z^{-1}\theta(z;p).
\end{equation}

Now we multiply the left-hand sides and right-hand sides of the first
identities in (\ref{eq-1}) and (\ref{eq-2}), respectively, and do the
same with the second identities. This yields the difference equations
\begin{equation}
\Gamma(qz;q,p)=\theta(z;p)\Gamma(z;q,p), \qquad
\Gamma(pz;q,p)=\theta(z;q)\Gamma(z;q,p),
\label{eq-3}\end{equation}
for the elliptic gamma function $\Gamma(z;q,p)$; in the explicit form,
we have
\begin{equation}
\Gamma(z;q,p) = \prod_{j,k=0}^\infty\frac{1-z^{-1}q^{j+1}p^{k+1}}{1-zq^jp^k}.
\label{ell-gamma}\end{equation}
Despite of the fact that the general idea of associating a generalized gamma
function with the elliptic theta function was formulated in the
well-known paper \cite{jac:basic}, it did not get much
attention. However, Jackson's double infinite product was used
explicitly in the mathematical physics literature on integrable
models of statistical mechanics starting with the Baxter's work
on the eight vertex model; see \cite{bax:partition}. The name ``elliptic gamma
function" for the product (\ref{ell-gamma}) was proposed by Ruijsenaars
in the recent paper \cite{rui:first}, where he reintroduced the function
$\Gamma(z;q,p)$ anew and started a systematic investigation of its
properties. A further detailed analysis of this function was performed
by Felder and Varchenko in \cite{fel-var:elliptic}.

In order to compare the elliptic gamma function with the Barnes multiple
gamma function, in (\ref{ell-gamma}) we put
$$
z=e^{2\pi i\frac{ u}{\omega_2}},\quad
q=e^{2\pi i\frac{ \omega_1}{\omega_2}},\quad
p=e^{2\pi i\frac{ \omega_3}{\omega_2}},
$$
where $\omega_i$ are some constants satisfying the same constraints
for the quasiperiods as in (\ref{m-gamma}).
Then it is not difficult to see that the set of zeros and poles of
$\Gamma(z;q,p)$, viewed as a meromorphic function of the
variable $u$, coincides with the set of zeros and poles of
the following combination of Barnes $\Gamma_3$-functions:
\begin{equation}
\frac{\Gamma_3(u;\omega_1,\omega_2,\omega_3)\Gamma_3(u-\omega_2;\omega_1,
-\omega_2,\omega_3)}{\Gamma_3(\omega_1+\omega_3-u;\omega_1,\omega_2,\omega_3)
\Gamma_3(\omega_1+\omega_3-\omega_2-u;\omega_1,-\omega_2,\omega_3)}.
\label{ell-g-bar1}\end{equation}
This means that the ratio of $\Gamma(z;q,p)$ and (\ref{ell-g-bar1}) is an
entire function of $u$, and this function is seen to be given by an
exponential of some polynomial of $u$ of the third degree.

For arbitrary complex $s$, the elliptic shifted factorials are defined
as ratios of elliptic gamma functions
$$
\theta(z;p;q)_s = \frac{\Gamma(zq^s;q,p)}{\Gamma(z;q,p)}.
$$
We use also the following shorthand notation:
\begin{eqnarray*}
&& \Gamma(t_1,\ldots,t_k;q,p)\equiv \prod_{j=1}^k \Gamma(t_j;q,p),
\\ &&
\theta(t_1,\ldots,t_k;p;q)_n\equiv\prod_{j=1}^k\prod_{\ell=0}^{n-1}
\theta(t_jq^{\ell};p), \quad n\in\mathbb{N}.
\end{eqnarray*}

The elliptic beta integral \cite{spi:elliptic} is the first exact
integration formula involving the elliptic gamma function.
We conclude this section by an explicit description of it.

\begin{theorem}
Let five complex parameters $t_m, m=0, \dots,4,$
satisfy the inequalities $|t_m|<1,\; |pq|<|A|,$
where $A\equiv\prod_{r=0}^4t_r.$
Define the elliptic beta integral as the following contour integral:
\begin{equation}\label{ell-int}
\mathcal{N}_E(\mathbf{t})=
\int_\mathbb{T}\Delta_E(z,\mathbf{t})\frac{dz}{z},
\end{equation}
where $\mathbb{T}$ is the positively oriented unit circle, and
\begin{equation}\label{weight}
\Delta_E(z,\mathbf{t}) = \frac{1}{2\pi i}
\frac{\prod_{m=0}^4\Gamma(zt_m, z^{-1}t_m; q,p)}
{\Gamma(z^2,z^{-2}, zA, z^{-1}A;q,p)}.
\end{equation}
Then
\begin{equation}\label{result}
\mathcal{N}_E(\mathbf{t})=\frac{2\prod_{0\leq m<s\leq 4} \Gamma(t_mt_s;q,p)}
{(q;q)_\infty(p;p)_\infty\prod_{m=0}^4\Gamma(At_m^{-1};q,p)}.
\end{equation}
\end{theorem}

Relation (\ref{ell-int}) determines a new Askey-Wilson type integral
representing an elliptic extension of the Nassrallah-Rahman integral.
Indeed, if we set $p=0$, then the integral (\ref{ell-int}) is
reduced to the Nassrallah-Rahman $q$-beta integral
\cite{nas-rah:projection,rah:integral} which, in turn, is a one
parameter extension of the celebrated Askey-Wilson $q$-beta integral
\cite{ask-wil:some}.
Theorem 1 was proved by the author with the help of an elliptic
generalization of the method used by Askey in \cite{ask:beta}
for proving the Nassrallah-Rahman integral. A large list of known
plain and $q$-hypergeometric beta integrals was given in
\cite{rah-sus:pearson}.

As shown in \cite{die-spi:elliptic}, a special finite-dimensional reduction
of (\ref{ell-int}) associated with the residue calculus results in the
elliptic generalization of the Jackson sum for a terminating $_8\Phi_7$
basic hypergeometric series, which was discovered by Frenkel and Turaev
in \cite{fre-tur:elliptic}. In \cite{spi:special}, the integral
(\ref{ell-int}) was applied to the construction of a large
family of continuous biorthogonal functions generalizing the Rahman's
$_{10}\Phi_9$ biorthogonal rational functions
\cite{rah:integral,rah:biorthogonality}.
These functions are expressed through products of two $_{12}E_{11}$
elliptic hypergeometric series with different modular parameters
(for the definition of an appropriate system of notations for such series,
see \cite{spi:theta}). It is believed that in the theory of biorthogonal
functions they play a role similar to that played by the Askey-Wilson
polynomials \cite{ask-wil:some} in the theory of orthogonal polynomials.
We describe these biorthogonal functions in the last section
and in Appendix A of the present paper and give complete proofs
of some results announced in \cite{spi:special}.
An elliptic extension of Wilson's discrete (finite-dimensional)
set of biorthogonal rational functions \cite{wil:orthogonal}
was constructed earlier by Zhedanov and the author in
\cite{spi-zhe:spectral,spi-zhe:classical}. The corresponding three
terms recurrence relation generates the most general example in the
pool of known terminating continued fractions expressed in terms of the
series of hypergeometric type, namely, it is expressed via the
$_{12}E_{11}$ series as well \cite{spi-zhe:gevp,spi-zhe:theory}.
The integral (\ref{ell-int}) leads to integral representations
for a terminating $_{12}E_{11}$ series and its particular bilinear
form \cite{spi:special} (the proofs are given in Appendix B).
All these results open new ways of exploration of the world
of elliptic functions and modular forms, which complement recent
progress reached in the classical setting by Milne \cite{mil:infinite}.

\bigskip
\centerline{\Large Acknowledgments}
\medskip

This paper is dedicated to Mizan Rahman, a master of $q$-special functions
from whom the author has learned much over ten years of regular contacts.
The author is also grateful to J.F. van Diejen for many useful discussions
during the work on \cite{die-spi:elliptic}-\cite{die-spi:review}, and for
remarks to this paper. Permanent encouragement from A.S. Zhedanov
is highly appreciated as well. Special thanks go
to C. Krattenthaler for explaining some of the Warnaar
theorems \cite{war:summation}, to V. Tarasov for emphasizing the
importance of the Barnes work \cite{bar:multiple} and for drawing
attention to the paper \cite{tar-var:geometry}, and to S. Kharchev
for a discussion of the properties of the generalized gamma functions.
Some of the main results were obtained during the author's stay at the
Max-Planck-Institut f\"ur Mathematik (Bonn) in the summer of 2002.
The author is indebted to this institute for warm hospitality and to
Yu.I. Manin, A. Okounkov, and D. Zagier for several stimulating
discussions at MPI of the subject of the present paper.
The organizers of the workshops ``Special Functions
in the Digital Age" (Minneapolis, USA, July 22-August 2, 2002) and
``Classical and Quantum Integrable Systems" (Protvino, Russia, January
9-11, 2003) are thanked for giving an opportunity to present there
some results of this work. This research is supported in part by the
Russian Foundation for Basic Research (RFBR) Grant No. 03-01-00781.

\section{A general definition of theta hypergeometric integrals}

The right-hand side of (\ref{ell-int}) belongs to a general
class of integrals related to the series of hypergeometric type built
from Jacobi theta functions. In accordance with the theory of general
theta hypergeometric series developed in \cite{spi:theta}, we give the
following definition.

\begin{definition}
Let $C$ denote a smooth Jordan curve on the complex plane.
Let $\Delta(y_1,\ldots,$ $y_n)$ be a meromorphic function
of its arguments $y_1,\ldots,y_n$. Consider the (multiple) integrals
\begin{equation}\label{definition}
I_n=\int_{C}dy_1\dots\int_Cdy_n \; \Delta(y_1,\ldots,y_n)
\end{equation}
and the ratios
\begin{equation}
h_\ell(\mathbf{y})=\frac{\Delta(y_1,\ldots,y_\ell+1,\ldots,y_n)
}{\Delta(y_1,\ldots,y_\ell,\ldots,y_n)}.
\label{hl}\end{equation}

Then the integrals $I_n$ are called:

1) the {\em plain hypergeometric integrals} if
\begin{equation}
h_\ell(\mathbf{y})=R_\ell(\mathbf{y}),
\label{plain-int}\end{equation}
are rational functions of $y_1,\ldots, y_n$ for all $\ell=1,\ldots, n;$

2) the {\em $q$-hypergeometric integrals} if
\begin{equation}
h_\ell(\mathbf{y})=R_\ell(q^\mathbf{y}),
\label{q-int}\end{equation}
are rational functions of $q^{y_1},\ldots, q^{y_n}$, $q\in\mathbb{C},$
for all $\ell=1,\ldots, n;$

3) the {\em elliptic hypergeometric integrals} if
for all $\ell=1,\ldots,n$ the ratios $h_\ell(\mathbf{y})$
are elliptic functions of the variables $y_1,\ldots,y_n$
with periods $\sigma^{-1}$ and $\tau\sigma^{-1}, \, \text{Im}(\tau)>0$;

4) the general {\em theta hypergeometric integrals} if $h_\ell(\mathbf{y})$
and $1/h_\ell(\mathbf{y})$ are meromorphic functions obeying
the double quasiperiodicity conditions
\begin{eqnarray}\nonumber
&& h_\ell(y_1,\ldots,y_k+\sigma^{-1},\ldots,y_n)
=e^{\sum_{j=1}^na_{\ell k}(j)y_j+b_{\ell k}} h_\ell(\mathbf{y}),
\\ && h_\ell(y_1,\ldots,y_k+\tau\sigma^{-1},\ldots,y_n)
=e^{\sum_{j=1}^n c_{\ell k}(j)y_j+d_{\ell k}} h_\ell(\mathbf{y}),
\label{quasi-per}\end{eqnarray}
with the quasiperiodicity factors similar to those
for the Weierstrass sigma function (which is related to
$\theta_1(u)$ in a simple way, see \cite{whi-wat:course}).
\end{definition}

If we assume that the variables $y_1,\ldots,y_n$ are discrete,
$\mathbf{y}\in\mathbb{N}^n$, and replace integrals by sums
$\sum_{\mathbf{y}\in\mathbb{N}^n}$, then we get the definitions
of the plain and $q$-hypergeometric series, which go back to
Horn \cite{ggr:general}, and the definition of the elliptic hypergeometric
series suggested in \cite{spi:theta}, respectively. The theta hypergeometric
series were defined in \cite{spi:theta} in a less general form
because of the less general choice of quasiperiodicity factors.
Evidently, if $a_{\ell k}(j)=b_{\ell k}=c_{\ell k}(j)=d_{\ell k}=0$,
then the theta hypergeometric functions are reduced to the elliptic ones.
The integrals (or series) defined in this way do not form an algebra
because, in general, sums of hypergeometric integrals do not fit the
taken definition.

The shifts $y_\ell\to y_\ell+1$ in (\ref{hl}) may be replaced by
translations by an arbitrary constant $y_\ell\to y_\ell +\omega_1,$
$\omega_1=const.$ However, we can replace $\omega_1$ by 1 after
an appropriate rescaling of $y_\ell$, which results in a simple
deformation of the contour $C$ in (\ref{definition}).

Consider the case of $n=1$ in detail. The general rational function of
$y$ can be represented in the form
$$
R(y)=\frac{\prod_{j=1}^n(1-a_j+y)\prod_{j=n+1}^r(a_j-1-y)}
{\prod_{j=1}^m(b_j-1-y)\prod_{j=m+1}^s(1-b_j+y)}\, x,
$$
where $n,r,m,s$ are arbitrary integers, $x$ is an arbitrary complex
constant, and $a_j, b_j$ describe the positions of the zeros and poles
of $R(y)$. The equation $\Delta(y+1)=R(y)\Delta(y)$
has the following general solution:
\begin{equation}
\Delta(y)=\frac{\prod_{j=1}^m\Gamma(b_j-y)
\prod_{j=1}^n\Gamma(1-a_j+y)}{\prod_{j=m+1}^s\Gamma(1-b_j+y)
\prod_{j=n+1}^r\Gamma(a_j-y)}\, x^{y}\varphi(y),
\label{mej}\end{equation}
where $\Gamma(y)$ is the standard gamma function and $\varphi(y)$
is an arbitrary periodic function, $\varphi(y+1)=\varphi(y)$.
If we set $\varphi(y)=1$, then for an appropriate choice of the
contour $C$ the integral $I_1$ (see (\ref{definition})) is none other
than the Meijer function \cite{erd:higher}. In this case we have no
natural additional tools for fixing an infinite dimensional
(functional) freedom contained in the solution $\Delta(y)$.

In the $q$-case, in a similar way we can write
$$
R(q^y)=\frac{\prod_{j=1}^n(1-t_jq^y)\prod_{j=n+1}^r(1-t_j^{-1}q^{-y})}
{\prod_{j=1}^m(1-w_j^{-1}q^{-y})\prod_{j=m+1}^s(1-w_jq^y)}\, x.
$$
For $0<|q|<1$, the general meromorphic solution of the equation
$\Delta(y+1)=R(q^y)\Delta(y)$ is
\begin{equation}
\Delta(y)= \frac{\prod_{j=n+1}^r (t_j^{-1}q^{1-y};q)_\infty
\prod_{j=m+1}^s(w_jq^y;q)_\infty}
{\prod_{j=1}^n(t_jq^y;q)_\infty
\prod_{j=1}^m(w_j^{-1}q^{1-y};q)_\infty}\, x^{y}\varphi(y),
\label{q-mej}\end{equation}
where, again, $\varphi(y)$ is an arbitrary periodic function,
$\varphi(y+1)=\varphi(y)$. In this case, for $\varphi(y)=1$
the integral $I_1$ describes a $q$-Meijer function, which was
investigated by Slater in \cite{sla:generalized}.

For $|q|>1$, the equation $\Delta(y+1)=R(q^y)\Delta(y)$ has the
following general solution
$$
\Delta(y)= \frac{\prod_{j=1}^n(t_jq^{y-1};q^{-1})_\infty
\prod_{j=1}^m(w_j^{-1}q^{-y};q^{-1})_\infty}{\prod_{j=n+1}^r
(t_j^{-1}q^{-y};q^{-1})_\infty \prod_{j=m+1}^s(w_jq^{y-1};q^{-1})_\infty}
\, x^{y}\varphi(y),
$$
that is we have an effective $q\to q^{-1}$ replacement and a reshuffling
of parameters in (\ref{q-mej}).

We remind the reader that $q=e^{2\pi i \sigma}$. The parameter $\sigma$
gives a second scale,
which may be used for generating a natural additional restriction
upon $\Delta(y)$. The function $q^y$ is periodic under the shift
$y\to y+\sigma^{-1}$, and (\ref{q-mej}) satisfies the equation
$\Delta(y+\sigma^{-1})/\Delta(y)=x^{1/\sigma}\varphi(y+\sigma^{-1})/
\varphi(y)$. We can fix $\varphi(y)$ by demanding that
$$
\varphi(y+\sigma^{-1})=\tilde R(e^{2\pi i y})\varphi(y),
$$
where $\tilde R$ is another rational function of its argument. In
accordance with the periodicity condition $\varphi(y+1)=\varphi(y)$,
we have
$$
\tilde R(e^{2\pi i y})=
\frac{\prod_{j=1}^{n'}(1-\tilde t_je^{-2\pi iy})
\prod_{j=n'+1}^{r'}(1-\tilde t_j^{-1}e^{2\pi iy})}
{\prod_{j=1}^{m'}(1-\tilde w_j^{-1}e^{2\pi iy})
\prod_{j=m'+1}^{s'}(1-\tilde w_je^{-2\pi iy})},
$$
where $\tilde t_j$ and $\tilde w_j$ are arbitrary new parameters.
Note that we cannot multiply the function $\tilde R$ by terms like
$\rho e^{2\pi i k y}, k\in\mathbb{Z}, \rho\in\mathbb{C},$ if they are
different from 1, because then the periodicity condition for $\varphi(y)$
will be broken. For $|q|<1$, the general meromorphic solution of the
difference equation for $\varphi(y)$ is as follows:
\begin{equation}
\varphi(y)= \frac{\prod_{j=m'+1}^{s'}(\tilde w_je^{-2\pi iy};
\tilde q)_\infty \prod_{j=n'+1}^{r'}(\tilde q\tilde t_j^{-1}e^{2\pi iy};
\tilde q)_\infty} {\prod_{j=1}^{m'} (\tilde q\tilde w_j^{-1}e^{2\pi iy};
\tilde q)_\infty \prod_{j=1}^{n'}(\tilde t_je^{-2\pi iy};
\tilde q)_\infty}\, \tilde \varphi(y),
\label{tq-mej}\end{equation}
where $\tilde q=e^{-2\pi i/\sigma}$ is the modular partner of $q$.
Indeed, for $\text{Im}(\sigma)>0$ we have $\text{Im}(\sigma^{-1})<0$,
and (\ref{tq-mej}) is well defined. The function $\tilde\varphi(y)$
in (\ref{tq-mej}) is an arbitrary elliptic function with periods
$1$ and $\sigma^{-1}$. It is characterized uniquely by the position of
its poles and zeros in the fundamental parallelogram of periods containing
$2k-1$ free parameters, where $k$ is the order of $\tilde \varphi(y)$.
Thus, the space of solutions is not too large: it becomes
finite-dimensional (in the sense of the number of free parameters).

Consider the regime $|q|=1$. Denoting $\sigma=\omega_1/\omega_2$
and assuming that $\text{Re}(\sigma)>0$, we introduce the variable
$u=y\omega_1$. Now it is possible to choose the
parameters $t_j, \tilde t_j$, etc in a special way, so that the
infinite products $(t_jq^y;q)_\infty$,
$(\tilde t_je^{-2\pi i y};\tilde q)_\infty$, etc in (\ref{q-mej}) and
(\ref{tq-mej}) combine into the double sine functions
$S(u+g_j;\omega_1,\omega_2)$ for some $g_j$, where
\begin{equation}
S(u;\omega_1,\omega_2) = \frac{(e^{2\pi i u/\omega_2}; q)_\infty}
{(e^{2\pi iu/\omega_1}\tilde q; \tilde q)_\infty},
\label{2d-sin}\end{equation}
is a well defined function for $|q|\to 1$. Indeed, it can be checked that
the zeros and poles of (\ref{2d-sin}) coincide with the zeros and poles
of the function $\Gamma_2(\omega_1+\omega_2-u;\mathbf{\omega})/\Gamma_2(u;
\mathbf{\omega})$, which is a well-defined meromorphic function of $u$
for $\omega_1/\omega_2>0$.

In this case $\sigma$ is real, and if it is incommensurate with 1, then
$\tilde\varphi(y)=1$ (i.e., the function $\Delta(y)$ is determined
quite uniquely). For a
description of the properties of the double sine function
and some of its applications, see \cite{jim-miw:quantum,kls:unitary,
nis-uen:integral,rui:generalized}. In particular,
the integrals introduced by Jimbo and Miwa in \cite{jim-miw:quantum}
as solutions of some $q$-difference equations at $|q|=1$ provided
the first examples of $q$-hypergeometric integrals for $q$ on the
unit circle. Faddeev's concept of the modular double for quantum groups
(see \cite{fad:modular}) is also related to the function (\ref{2d-sin}).

Thus, the world of $q$-Meijer functions appears to be reacher than in
the plain hypergeometric case. The introduction of the
additional equation involving shifts by $\sigma^{-1}$ brought some new
non-trivial structures in the integrals and reduced the functional
freedom in the definition of meromorphic function $\Delta(y)$ to an elliptic
function $\tilde\varphi(y)$ containing a finite number of free parameters.

Now, we turn to the single variable elliptic hypergeometric integrals.
The general elliptic function of order $r+1$ can be factorized as follows
\cite{whi-wat:course}:
\begin{eqnarray}
&& h(y)=e^\gamma\prod_{j=0}^r\frac{\theta_1(u_j+y;\sigma,\tau)}
{\theta_1(v_j+y;\sigma, \tau)}=e^\gamma
\frac{\theta(t_0q^y,\ldots,t_rq^y ;p)}
{\theta(w_0q^y,\ldots,w_rq^y;p)},
\label{n=1} \\
&& \makebox[4em]{}
\theta(t_0,\ldots, t_k;p)=\prod_{i=0}^k\theta(t_i;p),
\nonumber\end{eqnarray}
where $p=e^{2\pi i\tau }, \; \text{Im}(\tau)>0, \; q=e^{2\pi i\sigma}$.
The parameter $\gamma$ is an arbitrary complex number, but
$t_i\equiv q^{u_i}, \; w_i\equiv q^{v_i}$
satisfy the balancing constraint
\begin{equation}
\sum_{i=0}^r(u_i-v_i)=0, \quad \text{or}\quad
\prod_{i=0}^rt_i=\prod_{i=0}^rw_i,
\label{balance}\end{equation}
which guarantees that the meromorphic function $h(y)$ is doubly
periodic:
$$
h(y+\sigma^{-1})=h(y),\qquad   h(y+\tau\sigma^{-1})=h(y).
$$
For $\tau=\sigma$ (which requires that $\text{Im}(\sigma)> 0$), the
function $h(y)$ gives an explicit form of $\tilde\varphi(y)$ in
(\ref{tq-mej}).

In order to find the integrand $\Delta(y)$, it is necessary to solve
the first order difference equation
\begin{equation}
\Delta(y+1)=h(y)\Delta(y)
\label{1-eq}\end{equation}
in the class of meromorphic functions. The theory of such equations was
developed long ago (see, e.g., \cite{bar:linear}). Obviously, since
$h(y)$ is factorized into the ratio of products of theta functions,
it suffices to find a meromorphic solution of the equation
\begin{equation}
f(y+1)=\theta(q^y;p)f(y),
\label{ell-g-eq}\end{equation}
which leads to various elliptic gamma functions \cite{jac:basic}.
The simplest such function (\ref{ell-gamma}) is defined
from equation (\ref{ell-g-eq}) only up to a periodic function
$\varphi(y+1)=\varphi(y)$ and, moreover, it requires that
$\text{Im}(\sigma)>0$ (or $|q|<1$), which was not assumed in (\ref{n=1}).

We introduce the variable $z=q^y$, so that the shift $y\to y+1$ becomes
equivalent to the multiplication $z\to qz$. Then the general solution
of (\ref{1-eq}) looks like this:
\begin{equation}
\Delta(y)=\prod_{j=0}^r\frac{\Gamma(t_jz;q,p)}
{\Gamma(w_jz;q,p)} e^{\gamma y+\delta}\varphi(y),
\label{D}\end{equation}
where the balancing condition (\ref{balance}) is assumed
and $\varphi(y+1)=\varphi(y)$ is an arbitrary periodic function.
Using the reflection formulas
\begin{eqnarray}\nonumber
&& \Gamma(pz, qz^{-1};q,p)=\Gamma(qz,pz^{-1};q,p)
\\ && \makebox[2em]{} =\Gamma(pqz,z^{-1};q,p)=\Gamma(z,pqz^{-1};q,p)=1,
\label{refl-eq}\end{eqnarray}
in (\ref{D}) we can replace several elliptic gamma functions
containing in the arguments $z$ by those with arguments containing
$z^{-1}$. After that, $\Delta(y)$ would look closer to the integrands for
the plain or $q$-Meijer functions, but in the elliptic case this does
not increase generality because of the right-hand side of (\ref{refl-eq})
is trivial.

In the region $\text{Im}(\sigma)<0$, that is, for $|q|>1$,
the general solution of (\ref{1-eq}) can be written in the form
\begin{equation}
\Delta(y)=\prod_{j=0}^r\frac{\Gamma(w_jq^{y-1};q^{-1},p)}
{\Gamma(t_jq^{y-1};q^{-1},p)} e^{\gamma y+\delta}\varphi(y).
\label{D'}\end{equation}
Effectively, we have a permutation of parameters and a simple $q\to q^{-1}$
substitution in the elliptic gamma functions in (\ref{D}) (cf. the
definition of this function for $|q|>1$ given in \cite{fel-var:elliptic}).

Let us take $\varphi(y)=1$. Then the function (\ref{D}) satisfies
two simple difference equations of the first order:
\begin{eqnarray} \label{2-eq}
&&\Delta(y+\sigma^{-1})=e^{\gamma/\sigma}\Delta(y),  \\
&&\Delta(y+\tau\sigma^{-1})=e^{\gamma\tau/\sigma}
\prod_{j=0}^r \frac{\theta(t_jq^y;q)}{\theta(w_jq^y;q)}\Delta(y).
\label{3-eq}\end{eqnarray}
Suppose that $1, \sigma^{-1}, \tau\sigma^{-1}$ are pairwise incommensurate.
Then the system of three equations (\ref{1-eq}), (\ref{2-eq}), and
(\ref{3-eq}) determines $\Delta(y)$ uniquely up to a factor.
As in the $q$-hypergeometric case, we can generalize
equations (\ref{2-eq}) and (\ref{3-eq}), use them as natural tools
for fixing the functional freedom in $\Delta(y)$, and get
qualitatively different elliptic hypergeometric integrals in this way.

The ratio $\Delta(y+\tau\sigma^{-1})/\Delta(y)$ in (\ref{3-eq})
is an elliptic function with periods 1 and $\sigma^{-1}$.
Therefore, it is natural to demand that $\Delta(y+\sigma^{-1})/\Delta(y)$
be also an elliptic function with periods that, by symmetry, are
equal to 1 and $\tau\sigma^{-1}$.

\begin{theorem}
Suppose that $\Delta(y)$ satisfies equation (\ref{1-eq}) and
that 1, $\sigma^{-1}$, $\tau\sigma^{-1}$ are pairwise incommensurate.
Denote $\tilde q=e^{-2\pi i/\sigma}$, $\tilde p=e^{2\pi i\tau/\sigma}.$
For simplicity, assume that $\text{Im}(\sigma)>0$ (i.e., $|q|<1$).
If $\Delta(y+\tau\sigma^{-1})/\Delta(y)$ is an elliptic
function with periods 1 and $\tau\sigma^{-1}$, then for
$\text{Im}(\tau/\sigma)>0$ the most general form of the meromorphic
function $\Delta(y)$ is as follows:
\begin{equation}\label{delta-ell}
\Delta(y)=\prod_{j=0}^r\frac{\Gamma(t_jq^y;q,p)}{\Gamma(w_jq^y;q,p)}
\prod_{j=0}^{n}\frac{\Gamma(\tilde t_je^{-2\pi iy};\tilde q,\tilde p)}
{\Gamma(\tilde w_je^{-2\pi iy};\tilde q,\tilde p)}\,e^{\gamma y+\delta},
\end{equation}
where $\prod_{j=0}^rt_jw_j^{-1}=\prod_{j=0}^{n}\tilde t_j\tilde w_j^{-1}=1$.
For $\text{Im}(\tau/\sigma)<0$, we have
\begin{equation}\label{delta-ell'}
\Delta(y)=\prod_{j=0}^r\frac{\Gamma(t_jq^y;q,p)}{\Gamma(w_jq^y;q,p)}
\prod_{j=0}^{n}\frac{\Gamma(\tilde w_je^{-2\pi iy}\tilde p^{-1};
\tilde q,\tilde p^{-1})} {\Gamma(\tilde t_je^{-2\pi iy}\tilde p^{-1};
\tilde q,\tilde p^{-1})}\,e^{\gamma y+\delta}.
\end{equation}
\end{theorem}
\begin{proof}
First, observe that for $\text{Im}(\sigma)>0$ we have
$\text{Im}(\sigma^{-1})<0,$ automatically, that is $|\tilde q|<1.$
Therefore, for $\text{Im}(\tau/\sigma)>0$ the function
$\Gamma(z;\tilde q,\tilde p)$ is well defined.

The function $\Delta(y)$ in (\ref{D}) gives the general solution of
equation (\ref{1-eq}). Suppose that $\Delta(y+\sigma^{-1})/\Delta(y)$ is
an elliptic function of order $n+1$ with periods 1, $\tau\sigma^{-1}$.
For $\text{Im}(\tau/\sigma)>0$, this demand is equivalent to the
following equation for $\varphi(y)$:
\begin{equation}
\frac{\varphi(y+\sigma^{-1})}{\varphi(y)}=
\prod_{j=0}^{n}\frac{\theta(\tilde t_je^{-2\pi iy}; \tilde p)}
{\theta(\tilde w_je^{-2\pi iy}; \tilde p)},
\label{phi-eq}\end{equation}
where $\prod_{j=0}^{n}\tilde t_j\tilde w_j^{-1}=1$.
Note that we cannot multiply the right-hand side of this equation
by any constant different from 1, since this would violate the
condition $\varphi(y+1)=\varphi(y)$.

The meromorphic solution of (\ref{phi-eq}) is
\begin{equation}
\varphi(y)=\prod_{j=0}^{n}\frac{\Gamma(\tilde t_je^{-2\pi iy};\tilde q,
\tilde p)}{\Gamma(\tilde w_je^{-2\pi iy};\tilde q,\tilde p)}\,
\tilde \varphi(y),
\label{phi-form}\end{equation}
where $\tilde \varphi(y)$ is an elliptic function with periods 1 and
$\sigma^{-1}$. We can write
$$
\tilde\varphi(y)=\prod_{j=1}^m\frac{\theta(a_je^{-2\pi iy};\tilde q)}
{\theta(b_je^{-2\pi iy};\tilde q)}=\prod_{j=1}^m\frac{\Gamma(a_je^{-2\pi iy}
\tilde p,b_je^{-2\pi iy};\tilde q,\tilde p)}
{\Gamma(a_je^{-2\pi iy},b_je^{-2\pi iy}\tilde p;\tilde q,\tilde p)},
$$
where $\prod_{j=1}^ma_jb_j^{-1}=1$. Therefore, we can absorb the function
$\tilde\varphi(y)$ into the ratio of elliptic gamma functions in
(\ref{phi-form}) by changing $n\to n+2m$ and identifying
$\tilde t_k=\tilde p a_k, \tilde w_k=a_k$ for $k=n+1,\ldots, n+m$
and $\tilde t_k=\tilde b_k, \tilde w_k=\tilde p \tilde b_k$ for
$k=n+m+1,\ldots,n+2m$. Since $n$, $\tilde t_j$, $\tilde w_j$ are
arbitrary, without loss of generality we can set $\tilde\varphi(y)=1$,
which yields the desired expression (\ref{delta-ell}).

The function (\ref{delta-ell}) satisfies the following equations:
\begin{eqnarray} \label{2-eq'}
&&\Delta(y+\sigma^{-1})=e^{\gamma/\sigma}
\prod_{j=0}^{n}\frac{\theta(\tilde t_je^{-2\pi iy};\tilde p)}
{\theta(\tilde w_je^{-2\pi i y};\tilde p)}\, \Delta(y),  \\
&&\Delta(y+\tau\sigma^{-1})=e^{\gamma\tau/\sigma}
\prod_{j=0}^r \frac{\theta(t_jq^y;q)}{\theta(w_jq^y;q)}
\prod_{j=0}^{n}\frac{\theta(\tilde w_je^{-2\pi iy}\tilde p^{-1};\tilde q)}
{\theta(\tilde t_je^{-2\pi iy}\tilde p^{-1};\tilde q)}\, \Delta(y).
\label{3-eq'}\end{eqnarray}
The elliptic functions defined by the products $\prod_{j=0}^r$ and
$\prod_{j=0}^{n}$ in (\ref{3-eq'}) have different forms though both have
periods 1 and $\sigma^{-1}$. They are related to
each other by the modular transformation $\sigma\to -1/\sigma$
for the corresponding theta functions.

Now, we consider the region $\text{Im}(\tau/\sigma)<0$. Equations
(\ref{1-eq}) and (\ref{3-eq'}) are well defined in this case. They
can be used for the determination of $\Delta(y)$, and it can be
checked that, indeed, the function (\ref{delta-ell'}) provides
their general solution. Equation (\ref{2-eq'}) is replaced now
by the following one:
\begin{equation}
\frac{\Delta(y+\sigma^{-1})}{\Delta(y)}= e^{\gamma/\sigma} \prod_{j=0}^{n}
\frac{\theta(\tilde w_je^{-2\pi iy};\tilde p^{-1})}
{\theta(\tilde t_je^{-2\pi iy}; \tilde p^{-1})},
\label{phi-eq'}\end{equation}
that is, $\tilde p$ in (\ref{phi-eq}) is changed to $\tilde p^{-1}$,
and the parameters $\tilde t_j, \tilde w_j$ are replaced by
$\tilde p^{-1}\tilde w_j$ and $\tilde p^{-1}\tilde t_j$, respectively.
Using (\ref{D'}), it is easy to construct $\Delta(y)$ satisfying
(\ref{1-eq}), (\ref{2-eq'}), and (\ref{3-eq'}) in the $|q|>1$ region as well.
\end{proof}

In order to be able to work with $q$ on the unit circle $|q|=1$, we
need another elliptic gamma function: a kind of the elliptic analog of
the double sine function (\ref{2d-sin}). Denote
\begin{eqnarray}\nonumber
&& q=e^{2\pi i\frac{\omega_1}{\omega_2}}, \qquad
\tilde q =e^{-2\pi i\frac{\omega_2}{\omega_1}},
\\
&& p=e^{2\pi i\frac{\omega_3}{\omega_2}}, \qquad
\tilde p=e^{2\pi i\frac{\omega_3}{\omega_1}},
\label{ell-bases}\end{eqnarray}
where $\omega_i$ are some complex numbers.

Suppose that $\omega_1/\omega_2>0$ and $\text{Im}(\omega_3/\omega_2)>0$
(i.e., $|p|<1$). Then we have $\text{Im}(\omega_3/\omega_1)=
(\omega_2/\omega_1)\text{Im}(\omega_3/\omega_2)>0$, that is, $|\tilde p|<1$
automatically. Therefore, in the analysis of equation (\ref{1-eq})
for $|q|=1$ it is necessary to assume that $|p|, |\tilde p|<1$.

\begin{definition}
Let $|q|,|p|, |\tilde p|<1$. Then we define a new elliptic
gamma function by the formula
\begin{equation}
G(u;\mathbf{\omega})= \prod_{j,k=0}^\infty
\frac{(1-e^{-2\pi i\frac{u}{\omega_2}}q^{j+1}p^{k+1})
(1-e^{2\pi i \frac{u}{\omega_1}}{\tilde q}^{j+1}{\tilde p}^k)}
{(1-e^{2\pi i \frac{u}{\omega_2}}q^jp^k)
(1-e^{-2\pi i \frac{u}{\omega_1}}{\tilde q}^j{\tilde p}^{k+1})}.
\label{ell-d}\end{equation}
\end{definition}

In the limit $p\to 0$ taken in such a way that, simultaneously,
$\tilde p\to 0$, we get
$G(u;\omega_1,\omega_2,\omega_3)\to
S^{-1}(u;\omega_1,\omega_2)$, where the double sine function
$S(u;\mathbf{\omega})$ is fixed in (\ref{2d-sin}).

The function $G(u;\mathbf{\omega})$ satisfies
the following three difference equations:
\begin{eqnarray}
&& G(u+\omega_1;\mathbf{\omega})=\theta(e^{2\pi i\frac{u}{\omega_2}};p)
G(u;\mathbf{\omega}),
\label{ell-1eq} \\
&& G(u+\omega_2;\mathbf{\omega})=\theta(e^{2\pi i\frac{u}{\omega_1}};\tilde p)
G(u;\mathbf{\omega}),
\label{ell-2eq} \\
&& G(u+\omega_3;\mathbf{\omega})=
\frac{\theta(e^{2\pi i\frac{u}{\omega_2}};q)}
{\theta(e^{2\pi i\frac{u}{\omega_1}}\tilde q;\tilde q)}
G(u;\mathbf{\omega})
\nonumber \\ && \makebox[5.5em]{}
= S(u;\omega_1,\omega_2)
S(\omega_1+\omega_2-u;\omega_1,\omega_2)G(u;\mathbf{\omega}).
\label{ell-3eq}\end{eqnarray}
For pairwise incommensurate $\omega_1$, $\omega_2$, $\omega_3$,
these equations determine the meromorphic function $G(u;\mathbf{\omega})$
uniquely up to multiplication by a constant, which follows from
the nonexistence of triply periodic functions.

The first equation requires that $|p|<1$, the second one requires
$|\tilde p|<1$, and both of them do not impose any constraint
upon $q$. The third equation (\ref{ell-3eq}) involves only the
function $S(u;\omega_1,\omega_2)$, which is well defined for
$\omega_1/\omega_2>0$, that is, $|q|=|\tilde q|=1$. This means that
the function $G(u;\mathbf{\omega})$ may be well defined in this unit
circle region as well.

In essence, the original elliptic gamma function (\ref{ell-gamma})
has the same properties as the function (\ref{ell-g-bar1}).
In a similar way, the function (\ref{ell-d}) can be
expressed as the following combination of the Barnes
$\Gamma_3$-functions up to an exponential of some polynomial in $u$
of the third degree:
\begin{eqnarray}\nonumber
&& \frac{\Gamma_3(u;\omega_1,\omega_2,\omega_3)
\Gamma_3(\omega_3-\omega_1-u;-\omega_1,
-\omega_2,\omega_3)}{\Gamma_3(\omega_1+\omega_3-u;\omega_1,\omega_2,\omega_3)
\Gamma_3(u-\omega_1-\omega_2;-\omega_1,-\omega_2,\omega_3)}
\\ && \makebox[2em]{} \times
\Gamma_2(\omega_3-\omega_2-u;-\omega_2,\omega_3)
\Gamma_2(\omega_3-u;\omega_1,\omega_3).
\label{ell-g-bar2}\end{eqnarray}
From this representation it follows that, indeed, (\ref{ell-d}) is well
defined for real $\omega_1, \omega_2$ with $\omega_1/\omega_2>0$
(and any complex $\omega_3$), like in the
double sine function case. A more detailed analysis of this
correspondence and an investigation of other properties of the
function $G(u;\mathbf{\omega})$ will be given elsewhere.
In particular, it is expected that $G(u;\mathbf{\omega})$
is the key function for an elliptic extension of the modular
doubling principle for $q$-deformed algebras
\cite{fad:modular,kls:unitary}.

As a result, for $|q|=1$ we get a solution $\Delta(y)$ of equation
(\ref{1-eq}) by the mere replacement of $\Gamma(q^y;q,p)$ in (\ref{D})
by $G(y\omega_1;\mathbf{\omega})$. In the rest of this paper we limit
ourselves to the case where $|q|<1$. Note that the region $|p|=1$
is not well defined in the elliptic functions setting. In a sense, the
region of real $\omega_3/\omega_2$ is reachable only at the level of
the original Barnes multiple gamma functions.

\section{A theta analog of the Meijer function}

The integral corresponding to (\ref{D}) may be considered as a kind
of an elliptic extension of a particular Meijer function. The general
Jacobi theta functions analog of the Meijer function appears
in the case where $h(y)$ is a quasiperiodic function corresponding to the
fourth case of the definition given at the beginning of the previous
section, see (\ref{quasi-per}).

Let $P_3(y)=\sum_{i=1}^3\alpha_iy^i$ be an arbitrary polynomial of the
third degree obeying the property $P_3(0)=0$. The function defined by
the integral
\begin{equation}
G_r^s\left({\mathbf{t} \atop \mathbf{w}};\mathbf{\alpha};q,p\right)
=\int_C\frac{\prod_{j=0}^{s}\Gamma(t_jq^y;q,p)}
{\prod_{k=0}^r\Gamma(w_kq^y;q,p)} e^{P_3(y)}dy,
\label{e-meijer}\end{equation}
where $C$ is some contour on the complex plane, may be
called a theta analog of the Meijer function whenever the integral is
well defined. Note that no constraints are imposed in (\ref{e-meijer})
upon the integers $r, s$ and the complex parameters $t_j, w_k$

We have the following equation for the integrand $\Delta(y)$ of
(\ref{e-meijer}):
\begin{equation}
\frac{\Delta(y+1)}{\Delta(y)}=h(y)=e^{P_2(y)}
\frac{\theta(t_0q^y,\ldots,t_sq^y ;p)}
{\theta(w_0q^y,\ldots,w_rq^y;p)},
\label{theta-int}\end{equation}
where $P_2(y)=P_3(y+1)-P_3(y)$ is a polynomial in $y$ of the
second degree. From the considerations of \cite{spi:theta} it follows
that this $h(y)$ is the most general function such that $h$ is meromorphic
in $y$ (together with its inverse $1/h(y)$) and satisfies
the quasiperiodicity conditions
\begin{equation}
h(y+\sigma^{-1})=e^{ay+b}h(y),\qquad
h(y+\tau\sigma^{-1})=e^{cy+d} h(y)
\label{qua}\end{equation}
for some constants $a,b,c,d$. The function $h(y)$ may also be interpreted
as a general meromorphic modular Jacobi form in the sense of
Eichler and Zagier \cite{eic-zag:theory}.

However, the integral (\ref{e-meijer}) is not the most general integral
leading to (\ref{theta-int}). Using appropriate modifications of the
integrands (\ref{delta-ell}) and (\ref{delta-ell'}) and replacing $y$ by
$y/\omega_1$, we arrive at the general theta analog of the Meijer function.

\begin{definition}
In the definitions (\ref{ell-bases}) of the bases,
assume that $|q|, |p|<1$. Then, for $|\tilde p|<1$, the integral
\begin{eqnarray}\nonumber
\lefteqn{G_{rm}^{sn}\left({\mathbf{t},\mathbf{\tilde t}\atop \mathbf{w},
\mathbf{\tilde w}};\mathbf{\alpha};\mathbf{\omega}\right) } &&
\\ &&
=\int_C\frac{\prod_{k=0}^{s}\Gamma(t_ke^{\frac{2\pi iy}{\omega_2}};q,p)
\prod_{j=0}^{n}\Gamma(\tilde t_je^{-\frac{2\pi iy}{\omega_1}};
\tilde q,\tilde p)}{\prod_{k=0}^r\Gamma(w_ke^{\frac{2\pi iy}{\omega_2}};q,p)
\prod_{j=0}^m\Gamma(\tilde w_je^{-\frac{2\pi iy}{\omega_1}};
\tilde q,\tilde p)}\, e^{P_3(y)}dy
\label{gen-int}\end{eqnarray}
is called the general theta hypergeometric integral of one
variable whenever it is well defined. For $|\tilde p|>1$, we set
\begin{eqnarray}\nonumber
\lefteqn{G_{rm}^{sn}\left({\mathbf{t},\mathbf{\tilde t}\atop \mathbf{w},
\mathbf{\tilde w}};\mathbf{\alpha};\mathbf{\omega}\right) }&&
\\ &&
=\int_C\frac{\prod_{k=0}^{s}\Gamma(t_ke^{\frac{2\pi iy}{\omega_2}};q,p)
\prod_{j=0}^{m}\Gamma(\tilde w_je^{-\frac{2\pi iy}{\omega_1}}\tilde p^{-1};
\tilde q,\tilde p^{-1})}{\prod_{k=0}^r
\Gamma(w_ke^{\frac{2\pi iy}{\omega_2}};q,p)
\prod_{j=0}^n\Gamma(\tilde t_je^{-\frac{2\pi iy}{\omega_1}}\tilde p^{-1};
\tilde q,\tilde p^{-1})}\, e^{P_3(y)}dy.
\label{gen-int'}\end{eqnarray}
There are no constraints upon the integers $r, s, n, m \in\mathbb{N}$
and the complex parameters $t_j,\tilde t_j, w_k, \tilde w_k$.
\end{definition}

Both integrands of (\ref{gen-int}) and (\ref{gen-int'}) satisfy the
equations $\Delta(y+\omega_i)/\Delta(y)=h_i(y)$, $i=1,2,3,$ where $h_i$
are some quasiperiodic functions: $h_i(y+\omega_k)=e^{a_{ik} y+b_{ik}}h(y)$,
$i\neq k$, with $a_{ik}, b_{ik}$ being some constants related to the
parameters $\mathbf{t}, \mathbf{\tilde t},\mathbf{w},\mathbf{\tilde w},
\mathbf{\alpha}$ and $\mathbf{\omega}$. The integral (\ref{gen-int'})
was determined by the condition that it has the same functions $h_1(y),
h_3(y)$ as (\ref{gen-int}).
For a special choice of parameters $\mathbf{t}, \mathbf{\tilde t},\mathbf{w},
\mathbf{\tilde w},\mathbf{\alpha}$, in the limits $|p|, |\tilde p|\to 0$ or
$|p|, |\tilde p|^{-1}\to 0$ the function
$G_{rm}^{sn}\left({\mathbf{t},\mathbf{\tilde t}\atop \mathbf{w},
\mathbf{\tilde w}};\mathbf{\alpha};\mathbf{\omega}\right)$ is reduced
to the general $q$-hypergeometric integral considered in the previous
section, see (\ref{q-mej}) and (\ref{tq-mej}).

The general single variable theta hypergeometric series is defined
by the following formula \cite{spi:theta}:
\begin{eqnarray}\nonumber
\lefteqn{ {_{s+1}E_r}\left({t_0,\ldots, t_{s} \atop w_1,\ldots,w_r};
\mathbf{\alpha}; q,p\right) } &&
\\ && \makebox[4em]{}
= \sum_{n=0}^\infty \frac{\theta(t_0,t_1,\ldots,t_{s};p;q)_n}
{\theta(q,w_1,\ldots,w_r;p;q)_n}\, e^{P_3(n)}.
\label{_rE_s}\end{eqnarray}
Actually, these series are slightly more general than those
introduced in \cite{spi:theta}, because in that paper we
considered only the case where $\alpha_3=0$,
but the generalization to (\ref{_rE_s}) is straightforward.
We note that the presence of cubics of the independent variable
$y$ in (\ref{e-meijer}) or $n$ in (\ref{_rE_s}) is natural since we
are working at the level of the Barnes multiple gamma function
(\ref{m-gamma}) of the third order.

Writing (\ref{_rE_s}) in the form of the sum $\sum_{n=0}^\infty c_n$ with
$c_0=1$, we easily see that $c_{n+1}/c_n=h(n)$, where $h(n)$
is given by (\ref{theta-int}) with $w_0=q$ and $y=n$. This coincidence
is not artificial. Consider the sequence of poles of the integrand in
(\ref{e-meijer}) located at $y=y_0+n, n\in\mathbb{N}$, for some $y_0$.
We denote by $\kappa c_n, c_0=1,$ the residues of these poles. As
$y\to y_0+n$, we have $\Delta(y)\to \kappa c_n/(y-y_0-n) +O(1)$.
Now it is not difficult to see that
$$
\lim_{y\to y_0+n}\frac{\Delta(y+1)}{\Delta(y)}=\frac{c_{n+1}}{c_n}=
\lim_{y\to y_0+n} h(y)=h(y_0+n).
$$
In particular, this means that the sums of the residues in the integral
(\ref{e-meijer}) that appear from appropriate deformations of the
contour $C$, form the theta hypergeometric series (\ref{_rE_s}) for some
choices of the parameters.

In accordance with the classification introduced in \cite{spi:theta},
the elliptic hypergeometric series correspond, by definition, to $h(n)$
equal to an elliptic function of $n$. Such series are called also the
balanced theta hypergeometric series.
They are defined by the following constraints imposed upon (\ref{_rE_s}):
\begin{equation}
s=r, \qquad \alpha_3=\alpha_2=0, \qquad \prod_{j=0}^rt_j=\prod_{j=0}^r w_j.
\label{theta-balance}\end{equation}
Similarly, the integral (\ref{e-meijer}) will be called the {\em elliptic}
(or {\em balanced theta}) hypergeometric integral if the conditions
(\ref{theta-balance}) are satisfied. Evidently, in this case $h(y)$
in (\ref{theta-int}) becomes an elliptic function of $y$.

When $h(y)$ is an elliptic function of $y$ and of all parameters $u_j, v_j$
(we remind the reader that $t_j=q^{u_j}, w_j=q^{v_j}$), we call (\ref{_rE_s})
and (\ref{e-meijer}) the {\em totally elliptic hypergeometric series} and
integrals, respectively. As was shown in \cite{spi:theta},
in addition to the balancing requirement, such a property imposes the
following constraints on the parameters
\begin{equation}\label{well-poised-2}
t_jw_j=\rho=const, \qquad j=0,\ldots,r,
\end{equation}
which are known as the {\em well-poisedness} conditions in the theory of basic
hypergeometric series \cite{gas-rah:basic}. The explicit form of the
integrand $\Delta(y)$ for well-poised balanced theta hypergeometric
integrals is
\begin{equation}
\Delta(y)=\prod_{j=0}^{m-1}\frac{\Gamma(t_jz;q,p)}
{\Gamma(\rho t_j^{-1}z;q,p)}\; \frac{\Gamma(\rho^{\frac{m+1}{2}}
\prod_{j=0}^{m-1}t_j^{-1}\, z;q,p)}{\Gamma(\rho^{\frac{1-m}{2}}
\prod_{j=0}^{m-1}t_j\, z;q,p)}\,e^{\gamma y},
\label{wp-int}\end{equation}
where we have denoted $z=q^y$ and $\gamma=\alpha_1$.
The parameter $\rho$ is redundant, it can be eliminated by the rescalings
$t_i\to \rho^{1/2}t_i, z\to \rho^{-1/2}z$, but we keep it for further
needs. Observe that without loss of generality one of the parameters
in (\ref{e-meijer}) can be set equal to one by a shift of $y$.

Without the balancing condition, a theta hypergeometric series
$_{r+1}E_r$ is said to be {\em well-poised} if the constraints
(\ref{well-poised-2}) are valid with $w_0=q$, and
{\em very-well-poised} if, in addition to (\ref{well-poised-2}), we have
\begin{eqnarray}\nonumber
&& t_{r-3}=t_0^{1/2}q,\quad t_{r-2}=-t_0^{1/2}q, \\
&& t_{r-1}=t_0^{1/2}qp^{-1/2}, \quad t_{r} =- t_0^{1/2}qp^{1/2}.
\label{very-well-poised-2}\end{eqnarray}
Such series take a simpler form
\begin{eqnarray}\nonumber
\lefteqn{
_{r+1}E_r\left({t_0,t_1,\ldots, t_{r-4},qt_0^{1/2},-qt_0^{1/2},
qp^{-1/2}t_0^{1/2},-qp^{1/2}t_0^{1/2} \atop
qt_0/t_1,\ldots,qt_0/t_{r-4},t_0^{1/2},-t_0^{1/2},
p^{1/2}t_0^{1/2},-p^{-1/2}t_0^{1/2} };\mathbf{\alpha};q,p\right) } &&
\\ && \makebox[5em]{}
= \sum_{n=0}^\infty \frac{\theta(t_0q^{2n};p)}{\theta(t_0;p)}
\prod_{m=0}^{r-4}\frac{\theta(t_m;p;q)_n}{\theta(qt_0/t_m;p;q)_n}\, (-q)^n
e^{P_3(n)}.
\label{vwp-1}\end{eqnarray}
The essence of (\ref{very-well-poised-2}) consists in the replacement
of the product of four $\theta(t_iz;p)$ by one theta function
$\propto\theta(t_0q^2z^2;p)$ (this corresponds to doubling the
argument of the $\theta_1$-function). Very-well-poised series play a
distinguished role in applications, in particular, they admit an
appropriate generalization of the Bailey chains technique
of generating infinite
sequences of summation or transformation formulae \cite{spi:bailey}.

In the case of integrals, we call (\ref{e-meijer}) the
{\em very-well-poised} theta hypergeometric integral if, in addition to
conditions (\ref{well-poised-2}), eight parameters $t_i$ are fixed in
the following way:
\begin{eqnarray}
(t_{m-8},\ldots,t_{m-1})=(\pm(pq)^{1/2}, \pm q^{1/2}p, \pm p^{1/2}q,\pm pq).
\label{vwp}\end{eqnarray}
These constraints lead to squaring the argument of the
elliptic gamma function
\begin{equation}
\prod_{j=m-8}^{m-1}\Gamma(t_jz;q,p)=\Gamma(pqz^{2};q,p)=\frac{1}{\Gamma(z^{-2};q,p)}
\label{doubling}\end{equation}
(such a relation was used already in \cite{spi:elliptic} in the derivation
of the elliptic beta integral (\ref{ell-int})). As a result, the integrand
of the very-well-poised balanced theta hypergeometric integral takes the form
\begin{equation}
\Delta(y)=\prod_{j=0}^{m-9}\frac{\Gamma(t_jz;q,p)}
{\Gamma(\rho t_j^{-1}z;q,p)}\, \frac{\Gamma(\rho^{\frac{m+1}{2}}
p^{-6}q^{-6}\prod_{j=0}^{m-9}t_j^{-1}z;q,p)\,e^{\gamma y} }
{\Gamma(z^{-2},(\rho/pq)^2 z^2, \rho^{\frac{1-m}{2}}p^6q^6
\prod_{j=0}^{m-9}t_j\, z;q,p)}.
\label{vwp-int}\end{equation}
After imposing conditions (\ref{vwp}), the parameter $\rho$ is no longer
redundant, and its choice plays an important role. If we fix it as
$\rho=pq,$ then $\Delta(y)$ takes a particularly symmetric form
\begin{equation}
\Delta(y)=\frac{\prod_{j=0}^{m-9}\Gamma(t_jz,t_jz^{-1};q,p)}
{\Gamma(z^2,z^{-2}, Az, Az^{-1};q,p)}\,e^{\gamma y},
\label{vwp-int-fin}\end{equation}
where $A=(pq)^{\frac{13-m}{2}} \prod_{j=0}^{m-9}t_j$.
Clearly, the cases where $m$ is odd or even differ from each other in
a qualitative way. The choice $m=13$ gives the simplest expression
for $\Delta(y)$ and plays a distinguished role.
Other simple choices, $m=9$ or $11$, correspond to
particular subcases of the situation with $m=13$. For
$m=13$ and $\gamma =0$ we get the integrand $\Delta_E$ of
the elliptic beta integral (\ref{weight}), that is the
simplest very-well-poised elliptic hypergeometric integral
turns out to be exactly calculable when $C$ is taken to be a
special cycle corresponding to the unit circle on the $z$-plane.

The sums of the residues of the function (\ref{vwp-int-fin}) for $m=13$ are
expected to form a $_{14}E_{13}$ theta hypergeometric series.
However, the very-well-poisedness condition (\ref{vwp}) results
in the cancellation of theta functions in the corresponding ratios
of the series coefficients $h(n)$. As a result, we get only a $_{10}E_9$
very-well-poised elliptic hypergeometric series which, for
$\gamma=0$, corresponds to the left-hand side of (\ref{ft-sum}) or
(\ref{e-milne}) at $n=1$ (for more details, see
\cite{die-spi:elliptic,spi:theta}).
This shift of indices $m\to m-4$ brings in one more intriguing
point related to the origins of the very-well-poisedness condition. It
is necessary to find some deeper algebraic geometry explanations
of the fact that in the single variable case ``the nice things" (summation
or integration formulae) are related to the number 14, the order
of the initial elliptic function $h(y)$. As was shown in \cite{spi:theta,
spi:modularity}, in the multiple case this order raises to higher
numbers, but in quite an intriguing way as well.

Since $\Gamma((pq)^{1/2}z,(pq)^{1/2}z^{-1};q,p)=1$, we may drop
two parameters $\pm(pq)^{1/2}$ in the condition (\ref{vwp}). For
$\rho=pq, \gamma=0,$ this yields
\begin{equation}
\Delta(y)=\frac{\prod_{j=0}^{m-7}\Gamma(t_jz,t_jz^{-1};q,p)}
{\Gamma(z^2,z^{-2},-Az,-Az^{-1};q,p)},
\label{-A}\end{equation}
where $A=(pq)^\frac{11-m}{2} \prod_{j=0}^{m-7}t_j$.
For $m=11$ this expression looks similar to (\ref{weight}), but
the different sign in front of $A$ changes the things drastically,
and it is not known whether the corresponding integral gets any
closed form expression.

In a more general setting, we can impose balancing and very-well-poisedness
conditions upon general theta hypergeometric integrals (\ref{gen-int})
and (\ref{gen-int'}). However, at the moment it is not known whether
the simplest integrals appearing in this way admit exact evaluation.

As far as the multivariable integrals of hypergeometric type are
concerned, the general form of $\Delta(\mathbf{y})$ in the plain and
$q$-hypergeometric cases can be deduced from the Ore-Sato theorem
for Horn's series (see, e.g., \cite{ggr:general} for a detailed
discussion). The general form of the multiple elliptic hypergeometric
series or integrals is not established yet. We formulate it as an open
problem---to find an elliptic or general theta functions analog of
the Ore-Sato characterization theorem. In the following sections we
give a series of examples of multivariable extensions of the
very-well-poised balanced theta hypergeometric integrals
associated with the root systems $A_n$ and $C_n$.

As to the further possible
generalizations, it is natural to consider integrals of hypergeometric
type for arbitrary algebraic curves or general Abelian varieties.
Both would involve Riemann theta functions of many variables,
appropriate generalizations of gamma functions and theta hypergeometric
series. Some preliminary discussion of ideas in this direction can be
found in \cite{spi:modularity}, in particular, a special subcase
of the $_8\Phi_7$ Jackson summation formula was generalized
there to Riemann surfaces of arbitrary genus.

\section{Known $C_n$ elliptic beta integrals}

The following multivariable generalization of the Euler beta
integral (\ref{euler}) has been introduced by Selberg \cite{aar:special}:
\begin{eqnarray}
\lefteqn{\int_0^1\cdots\int_0^1
\prod_{1\leq j\leq n} x_j^{\alpha-1} (1-x_j)^{\beta-1}
\prod_{1\leq j<k\leq n}|x_j-x_k|^{2\gamma}\; dx_1\cdots dx_n } &&
\nonumber \\
&& =\prod_{1\leq j\leq n} \frac{\Gamma(\alpha+(j-1)\gamma)
\Gamma(\beta+(j-1)\gamma) \Gamma (1+j\gamma)}
     {\Gamma (\alpha+\beta+(n+j-2)\gamma) \Gamma (1+\gamma)},
\label{Sint}\end{eqnarray}
where $\text{Re}(\alpha ),\, \text{Re} (\beta) >0$, and
$\text{Re}(\gamma ) > - \min (1/n,
\text{Re} (\alpha )/(n-1) ,\text{Re} (\beta)/(n-1) )$.
This integral has found many important applications in
mathematical physics.

At the first glance, the integrand of (\ref{Sint})
does not fit the definition of the hypergeometric integrals introduced
in the preceding section. However, we can rescale $y_i\to y_i\epsilon $ in
(\ref{definition}), and choose rational functions $R_i(\mathbf{y})$
appropriately so that the limit $\epsilon\to 0$ becomes well
defined and yields a system of linear differential equations of
the first order. As a result, we get $\partial\Delta(\mathbf{y})/
\partial y_i=R_i(\mathbf{y})\Delta(\mathbf{y})$ and, apparently,
(\ref{Sint}) satisfies these conditions.

Two types of multidimensional generalizations of the elliptic
beta integral (\ref{ell-int}) to the root system $C_n$ were
proposed by van Diejen and the author in \cite{die-spi:elliptic,
die-spi:selberg}. One of them reduces in a special limit to
the Selberg integral (\ref{Sint}). Here we describe these
elliptic Selberg integrals explicitly. For brevity, we drop the
bases $p,q$ in the notation for elliptic gamma functions from now on.
We introduce the Type I integrand as
\begin{eqnarray}
\Delta^I ({\bf z};C_n)&=& \frac{1}{(2\pi i)^n}
\prod_{1\leq j<k\leq n}
\Gamma^{-1}(z_jz_k,z_jz_k^{-1},z_j^{-1}z_k,z_j^{-1}z_k^{-1}) \nonumber  \\
&&\times\prod_{j=1}^n\frac{\prod_{r=0}^{2n+2}\Gamma(t_rz_j,t_rz_j^{-1})}
{\Gamma(z_j^2,z_j^{-2}, Az_j,Az_j^{-1})},
\label{esintA}\end{eqnarray}
where $t_r\in\mathbb{C},\, r=0,\ldots,2n+2,$ are free parameters
and $A\equiv\prod_{r=0}^{2n+2}t_r$. The Type II integrand has the form
\begin{eqnarray}
\Delta^{II} (\mathbf{z};C_n)&=&\frac{1}{(2\pi i)^n}\prod_{1\leq j<k\leq n}
\frac{\Gamma(tz_jz_k,tz_jz_k^{-1},tz_j^{-1}z_k,tz_j^{-1}z_k^{-1})}
{\Gamma(z_jz_k,z_jz_k^{-1},z_j^{-1}z_k,z_j^{-1}z_k^{-1})} \nonumber \\
&& \times \prod_{j=1}^n\frac{\prod_{r=0}^4\Gamma(t_rz_j,t_rz_j^{-1})}
{\Gamma(z_j^2,z_j^{-2}, B z_j, B z_j^{-1})},
\label{esintB}\end{eqnarray}
where $t, t_r\in \mathbb{C},\, r=0,\ldots,4,$ are free parameters and
$B\equiv t^{2n-2}\prod_{s=0}^4t_s$. By $\mathbb{T}$ we denote the unit
circle with positive orientation.

Consider the first type of the $C_n$ multivariable elliptic beta integral.
We take $|p|,|q|<1$ and $|t_r|<1,\, r=0,\ldots, 2n+2,$ and
assume that $|pq|<|A|$. Then
\begin{eqnarray}
\lefteqn{\int_{\mathbb{T}^n}\Delta^I ({\bf z};C_n)\frac{dz_1}{z_1}
\cdots \frac{dz_n}{z_n} } &&  \nonumber \\ &&
=\frac{2^n n!}{(p;p)_\infty^n(q;q)_\infty^n}
\frac{\prod_{0\leq r<s\leq 2n+2}\Gamma(t_rt_s)}
{\prod_{r=0}^{2n+2}\Gamma(t_r^{-1}A)}.
\label{SintA}\end{eqnarray}

The second type of the multiple elliptic beta integral has the following
form. We take $|p|,|q|, |t|<1$ and $|t_r|<1,\, r=0,\ldots, 4,$ and
assume that $|pq|<|B|$. Then
\begin{eqnarray}
\lefteqn{\int_{\mathbb{T}^n}\Delta^{II} ({\bf z};C_n)\frac{dz_1}{z_1}
\cdots \frac{dz_n}{z_n} } && \nonumber \\
&& = \frac{2^n n!}{(p;p)_\infty^n(q;q)_\infty^n}
\prod_{j=1}^n\frac{\Gamma(t^j)}{\Gamma(t)}
\frac{\prod_{0\leq r<s\leq 4}\Gamma(t^{j-1}t_rt_s)}
{\prod_{r=0}^4\Gamma(t^{1-j}t_r^{-1}B)}.
\label{SintB}\end{eqnarray}

The elliptic Selberg integral of Type II (\ref{SintB}) can be deduced
from the Type I integral with the help of some interesting trick
\cite{die-spi:selberg}. The type I integral (\ref{esintA}) was
proved in \cite{die-spi:selberg} under a vanishing hypothesis,
namely, that its left-hand side vanishes on the hypersurface of parameters
$A=t_{2n+2}$ after an appropriate deformation of $\mathbb{T}$ to an
integration contour $C$ that separates the sequences of poles in
$\mathbf{z}$ converging to zero from those diverging to infinity.
As $p\to 0$, both Type I and II integrals are reduced to Gustafson's
well-known $q$-Selberg integrals \cite{gus:some1,gus:some2},
which are related (for particular choices of parameters) to the
Macdonald-Morris constant term identities and the Macdonald polynomials
for various root systems \cite{mac:constant}, including the
Koornwinder polynomials \cite{kor:pol}.

As was shown in \cite{die-spi:elliptic,
die-spi:modular}, the sums of residues of the functions (\ref{esintA})
and (\ref{esintB}) form some multiple elliptic hypergeometric series.
In particular, the multivariable $_{10}E_9$ sum conjectured by Warnaar
in \cite{war:summation} can be deduced from (\ref{SintB}).
The residue calculus for (\ref{SintA}) yields an elliptic extension
of the basic hypergeometric series summation formula proved
in \cite{den-gus:beta} and \cite{mil-lil:consequences}.
Recursive proofs of these multivariable Frenkel-Turaev sums
were given by Rosengren in \cite{ros:proof,ros:elliptic}. The property
of well-poisedness (or total ellipticity, in the case of elliptic
hypergeometric series) plays an important role in such summation
formulas. First examples of (plain) multiple hypergeometric series
well-poised on classical groups were introduced
by Biedenharn, Holman, and Louck in \cite{hbl:hypergeometric}.

\section{A new $C_n$ integration formula}

We define
\begin{eqnarray}\nonumber
&& \Delta^{III}(\mathbf{z};C_n) = \frac{1}{(2\pi i)^n}
\prod_{1\leq i<j\leq n}z_j\theta(z_iz_j^{-1},z_i^{-1}z_j^{-1};p)
\\ && \makebox[2em]{} \times
\prod_{i=1}^n\prod_{\nu=\pm 1}\frac{\Gamma(z_i^\nu x_i,z_i^\nu t_1,
z_i^\nu t_2,z_i^\nu t_3,z_i^\nu t/x_i)}
{\Gamma(z_i^{2\nu},z_i^\nu A)},
\label{delta-c3}\end{eqnarray}
where $A=tt_1t_2t_3q^{n-1}$.

\begin{theorem}
Impose the following restrictions upon parameters: $|x_i|,$ $|t_k|<1,$
$|t|<|x_i|$, where $i=1,\ldots, n$, $k=1,2,3$, and $|pq|<|A|$. Then
\begin{eqnarray}\nonumber
&& \int_{\mathbb{T}^n}\Delta^{III}(\mathbf{z};C_n)
\frac{dz_1}{z_1}\cdots\frac{dz_n}{z_n}
=\frac{2^n}{(p;p)_\infty^n(q;q)_\infty^n}
\prod_{1\leq i<j\leq n}x_j\theta\left(x_i/x_j,t/x_ix_j;p\right)
\\ && \makebox[4em]{} \times
\Gamma^n(t)\prod_{i=1}^n\left(\frac{\prod_{1\leq r<s\leq 3}
\Gamma(t_rt_sq^{i-1})}{\Gamma(A/x_i,Ax_i/t)} \prod_{k=1}^3\frac{
\Gamma(x_it_k,tt_k/x_i)}{\Gamma(Aq^{1-i}/t_k)}\right).
\label{c3}\end{eqnarray}
\end{theorem}

\begin{proof}
Consider the determinant
\begin{eqnarray}\nonumber
&& \det_{1\leq i,j\leq n}\left(\int_\mathbb{T}
\Delta_E(z,x_i,t_1q^{n-j},t_2q^{j-1},t_3,tx_i^{-1})\frac{dz}{z}\right)
\\  && \makebox[2em]{}
=\frac{1}{(2\pi i)^n}\int_{\mathbb{T}^n}\frac{dz_1}{z_1}\cdots
\frac{dz_n}{z_n}\prod_{i=1}^n G_i(z_i)G_i(z_i^{-1})\; D(\mathbf{z}),
\label{det1}\end{eqnarray}
where $\Delta_E$ is the integrand of the elliptic beta integral
(\ref{weight}) with an appropriate choice of the parameters. The expression
standing on the right-hand side of (\ref{det1}) appears after taking
the integral signs outside of the determinant symbol, so that we get
a multiple integral with
\begin{eqnarray*}
&& G_i(z_i)=\frac{\Gamma(z_ix_i,z_it_1,z_it_2,z_it_3,z_it/x_i)}
{\Gamma(z_i^2,z_iA)},
\\
&& D(\mathbf{z})=\det_{1\leq i,j\leq n}\left(\theta(z_it_1,
z_i^{-1}t_1;p;q)_{n-j}\theta(z_it_2,z_i^{-1}t_2;p;q)_{j-1}\right).
\end{eqnarray*}
The determinant $D(\mathbf{z})$ can be rewritten as follows
\begin{eqnarray*}
&& D(\mathbf{z})=\frac{\prod_{i=1}^n \theta(z_it_2,z_i^{-1}t_2;p;q)_{n-1}}
{t_2^{2\binom{n}{2}} q^{4\binom{n}{3}} }
\\ && \makebox[4em]{} \times
\det_{1\leq i,j \leq n}
\left(\frac{\theta(z_it_1,z_i^{-1}t_1;p;q)_{n-j}}
{\theta(q^{2-n}/z_it_2,q^{2-n}z_i/t_2;p;q)_{n-j}}\right).
\end{eqnarray*}
In \cite{war:summation}, Warnaar computed the following
elliptic generalization of a Krattenthaler determinant \cite{kra:major}:
\begin{eqnarray}
&& \det_{1 \leq i,j \leq n} \left(
\frac{ \theta(aX_i,ac/X_i;p;q)_{n-j} }
{ \theta(bX_i,bc/X_i;p;q)_{n-j} } \right)
\nonumber \\ && \makebox[4em]{}
= a^{\binom{n}{2}}q^{\binom{n}{3}}
\prod_{1\leq i<j\leq n}X_j\theta(X_iX_j^{-1},cX_i^{-1}X_j^{-1};p)
\nonumber \\ && \makebox[6em]{}
\times \prod_{i=1}^n \frac{\theta(b/a,abcq^{2n-2i};p;q)_{i-1}}
{\theta(bX_i,bc/X_i;p;q)_{n-1}}.
\label{e-kratt} \end{eqnarray}
Using this identity for $X_i=z_i$, $a=t_1, b=q^{2-n}/t_2$, and $c=1$,
we find
\begin{eqnarray}\nonumber
&& D(\mathbf{z})=(t_1t_2^2)^{\binom{n}{2}}q^{3\binom{n}{3}}
\prod_{1\leq i<j\leq n}z_j\theta(z_iz_j^{-1},z_i^{-1}z_j^{-1};p)
\\ && \makebox[4em]{} \times
\prod_{i=1}^n \theta(q^{2-n}/t_1t_2,t_1q^{n+2-2i}/t_2;p;q)_{i-1}.
\label{det-aux}\end{eqnarray}
As a result, the determinant (\ref{det1}) yields an expression
proportional to the left-hand side of the $C_n$ multiple integral in
question (\ref{c3}).

Now we substitute the result of computation of the elliptic beta integral
(\ref{ell-int}) into the determinant (\ref{det1}). This yields
\begin{eqnarray}\nonumber
&& \text{l.h.s. of (\ref{det1})} = \frac{2^n}{(q;q)_\infty^n(p;p)_\infty^n}
\\ && \makebox[1em]{}  \times
\prod_{i=1}^n\frac{\Gamma(x_it_1,t_1t/x_i,x_it_2,t_2t/x_i,
x_it_3,t,t_1t_2q^{n-1},t_1t_3q^{n-i},t_2t_3q^{i-1}, t_3t/x_i)}
{\Gamma(A/x_i,A/t_3,x_iA/t,Aq^{1-i}/t_2,Aq^{i-n}/t_1)}
\nonumber \\ && \makebox[1em]{}
\times \det_{1\leq i,j\leq n}\left(\theta(x_it_1,t_1t/x_i;p;q)_{n-j}
\theta(x_it_2,t_2t/x_i;p;q)_{j-1}\right).
\label{rhs}\end{eqnarray}

By (\ref{e-kratt}), the determinant in the
last line takes the form
\begin{eqnarray*}\nonumber
&& (t_1t_2^2t)^{\binom{n}{2}}q^{3\binom{n}{3}}
\prod_{1\leq i<j\leq n}x_j\theta\left(x_i/x_j,t/x_ix_j;p
\right) \\ && \makebox[4em]{} \times
\prod_{i=1}^n \theta(q^{2-n}/t_1t_2t,t_1q^{n+2-2i}/t_2;p;q)_{i-1}.
\end{eqnarray*}
Equating the resulting expression in (\ref{rhs}) with the right-hand side
of (\ref{det1}), we get the desired integral evaluation (\ref{c3}).
\end{proof}

This is the first nontrivial multiple elliptic beta integral with
a complete proof. It is not symmetric in $p$ and $q$, unlike
all other cases considered in the present paper. This fact suggests that
there should exist yet another
integral of similar nature that would be symmetric in $p,q$.

The above method of computation of the taken $C_n$ integral represents
a next step in the logical development of applications of determinant
formulas to multiple basic hypergeometric series; see, e.g., the paper
\cite{gus-kra:determinant} of Gustafson and Krattenthaler, which was
followed by Schlosser \cite{sch:summation,sch:nonterminating} and
Warnaar \cite{war:summation}. Similar considerations for computing
some multiple $q$-hypergeometric integrals were given by
Tarasov and Varchenko in \cite{tar-var:geometry}.

\section{An elliptic beta integral for the $A_n$ root system}

In this section we conjecture a multiple elliptic beta integral
for the $A_n$ root system, which will be used in the next section
for derivation of other nontrivial $A_n$ integrals.

\begin{conjecture}
Let $z_i,\: i=1,\ldots,n,$  $t_k,\: k=1,\ldots,$ $n+1,$ and
$f_j,\: j=1,\ldots,$ $n+2,$ be independent complex variables
($n$ is an arbitrary positive integer). We denote
$A\equiv \prod_{k=1}^{n+1}t_k$, $B\equiv\prod_{j=1}^{n+2}f_j$, and
\begin{equation}
\Delta^I(\mathbf{z};A_n)=\frac{1}{(2\pi i)^n}
\frac{\prod_{k=1}^{n+1}\left(\prod_{i=1}^{n+1}\Gamma(t_iz_k^{-1})
\prod_{j=1}^{n+2}\Gamma(f_jz_k)\right)}
{\prod_{i,j=1;\: i\neq j}^{n+1}\Gamma(z_iz_j^{-1})
\prod_{k=1}^{n+1}\Gamma(ABz_k)},
\label{delta-an}\end{equation}
where $z_1z_2\cdots z_{n+1}=1$.

Suppose that the parameters $t_k, f_j$ satisfy the constraints
$|t_k|, |f_j|<1, |pq|<|AB|$.
Then the following integration formula is conjectured to hold true:
\begin{eqnarray}\nonumber
&& \int_{\mathbb{T}^n}\Delta^I(\mathbf{z};A_n)
\frac{dz_1}{z_1}\cdots\frac{dz_n}{z_n}
=\frac{(n+1)!}{(q;q)_\infty^n(p;p)_\infty^n}
\\ && \makebox[2em]{} \times
\frac{\Gamma(A)\prod_{j=1}^{n+2}\Gamma(f_j^{-1}B)\:
\prod_{k=1}^{n+1}\prod_{j=1}^{n+2}\Gamma(t_kf_j)}
{\prod_{k=1}^{n+1}\Gamma(t_kB)\:\prod_{j=1}^{n+2}\Gamma(f_j^{-1}AB)}.
\label{int-an}  \end{eqnarray}
\end{conjecture}

For $n=1$ this conjecture is reduced to the elliptic beta integral
(\ref{ell-int}). For arbitrary $n$ and $p=0$, we get a Gustafson
integral proved in \cite{gus:some1}. Let us show that the two sides
of (\ref{int-an}) satisfy one and the same difference equation.

\begin{theorem}
Let $I_n(\mathbf{t},\mathbf{f})$ denote either side of (\ref{int-an}).
Then this function satisfies the $q$-difference equation
\begin{equation}
\sum_{r=1}^{n+1}
\frac{\theta(Bt_r;p)}{\theta(A;p)}
\prod_{\stackrel{j=1}{j \neq r}}^{n+1}
\frac{\theta(ABt_j;p)}{\theta(t_rt_j^{-1};p)}
I_n(t_1,\ldots,qt_r,\ldots,t_{n+1},\mathbf{f})=I_n(\mathbf{t},\mathbf{f})
\label{an-eq}\end{equation}
and its partner obtained by the permutation of $q$ and $p$.
\end{theorem}
\begin{proof}
Denote the function (\ref{delta-an}) by
$\Delta^I(\mathbf{z};t_1,\ldots,t_{n+1};A_n)$.
It is not difficult to see that
$$
\frac{\Delta^I(\mathbf{z};t_1,\ldots,qt_r,\ldots,t_{n+1};A_n)}
{\Delta^I(\mathbf{z};t_1,\ldots,t_{n+1};A_n)}
=\prod_{k=1}^{n+1}\frac{\theta(t_rz_k^{-1};p)}{\theta(ABz_k;p)},
$$
so that equation (\ref{an-eq}) for the left-hand side of (\ref{int-an})
is satisfied if the following theta functions identity is fulfilled:
\begin{equation}
\sum_{r=1}^{n+1} \frac{\theta(Bt_r;p)}{\theta(A;p)}
\prod_{\stackrel{j=1}{j \neq r}}^{n+1}
\frac{\theta(ABt_j;p)}{\theta(t_rt_j^{-1};p)}
\prod_{k=1}^{n+1}\frac{\theta(t_rz_k^{-1};p)}{\theta(ABz_k;p)}=1.
\label{id1}\end{equation}
For $n=1$ this identity is equivalent to the well-known relation
for products of four theta functions
\begin{eqnarray}\nonumber
\lefteqn{\theta(xw,x/w,yz,y/z;p) -\theta(xz,x/z,yw,y/w;p)}
\makebox[8em]{}  && \\
&& =yw^{-1}\theta(xy,x/y,wz,w/z;p)
\label{ident}\end{eqnarray}
and equation (\ref{an-eq}) coincides with that used in
\cite{spi:elliptic} for proving the integral (\ref{ell-int}).

For $n>1$, identity (\ref{id1}) can be established with the
help of the Liouville theorem, much as in the arguments presented in
\cite{die-spi:selberg}. A simpler proof follows from the general theta
functions identity given in \cite{whi-wat:course}.
As was shown by Rosengren in \cite{ros:elliptic}, that identity
can be rewritten as the following generalized partial fractions
expansion of a ratio of theta functions:
\begin{equation}
\prod_{k=1}^n\frac{\theta(t/b_k;p)}{\theta(t/a_k;p)}
=\sum_{r=1}^n\frac{\theta(t a_1\cdots a_n /a_rb_1\cdots b_n;p)}
{\theta(t/a_r,a_1\cdots a_n/b_1\cdots b_n;p)}
\frac{\prod_{j=1}^n\theta(a_r/b_j;p)}
{\prod_{\stackrel{j=1}{j \neq r}}^n\theta(a_r/a_j;p)},
\label{id2}\end{equation}
where $a_1\cdots a_n\neq b_1\cdots b_n$.
Here we replace $n$ by $n+1$ and substitute $a_k=t_k^{-1}$,
$b_k=z_k^{-1}$, and $t=AB$. As a result, we get an identity
which is seen to coincide with (\ref{id1}) due to the relation
$z_1\cdots z_{n+1}=1$.

In a similar way, for the right-hand side of (\ref{int-an}) we get
$$
\frac{I_n(t_1,\ldots,qt_r,\ldots,t_{n+1}, \mathbf{f} ) }
{I_n(\mathbf{t},\mathbf{f})}=
\frac{\theta(A;p)}{\theta(t_rB;p)}
\prod_{j=1}^{n+2}\frac{\theta(t_rf_j;p)}{\theta(ABf_j^{-1};p)},
$$
and in this case equation (\ref{an-eq}) becomes equivalent to the identity
\begin{equation}
\frac{\prod_{j=1}^{n+2}\theta(ABf_j^{-1};p)}
{\prod_{j=1}^{n+1}\theta(ABt_j;p)}= \sum_{r=1}^{n+1}
\frac{ \prod_{j=1}^{n+2}\theta(t_rf_j;p) }
{\prod_{\stackrel{j=1}{j \neq r}}^{n+1} \theta(t_rt_j^{-1};p) }
\frac{1}{\theta(ABt_r;p)}.
\label{id3}\end{equation}
If we substitute here $f_{n+2}=B/f_1\cdots f_{n+1}$ and divide both
sides by $\theta(ABf_{n+2}^{-1};p)$, then we get (\ref{id2})
with $n$ replaced by $n+1$ and with $a_j=t_j^{-1}, b_j=f_j, t=AB.$
\end{proof}

The equation derived above works in the space
of parameters $t_k$, whereas in \cite{gus:some1} Gustafson used
an equation in the variables $f_j$ for proving the $p=0$ case of
the integral (\ref{int-an}). It is of interest that the latter equation
does not admit a straightforward elliptic generalization, namely,
the corresponding partial fraction expansion cannot be lifted
to the theta functions level.

Another argument in favor of the validity of formula
(\ref{int-an}) consists in the fact that, via the residue calculus,
it generates a multivariable $_{10}E_9$
elliptic hypergeometric series sum for the $A_n$ root system, which
was proved by Rosengren in \cite{ros:elliptic} and which was
considered independently by the author in \cite{spi:modularity}.

\begin{theorem}
The residue calculus for the integral (\ref{int-an}) yields the
following summation formula:
\begin{eqnarray} \nonumber
\lefteqn{ \sum_{\stackrel{0\leq \lambda_j \leq N_j}{j=1,\ldots, n}}
q^{\sum_{j=1}^nj\lambda_j}
\prod_{j=1}^n\frac{\theta(t_jq^{\lambda_j+|\lambda|};p)}{\theta(t_j;p)}
\prod_{1\leq i<j \leq n} \frac{\theta(t_it_j^{-1}q^{\lambda_i-\lambda_j};p)}
{\theta(t_it_j^{-1};p)}  } &&
\\ \nonumber
&& \times \prod_{i,j=1}^n\frac{\theta(t_it_j^{-1}q^{-N_j};p;q)_{\lambda_i}}
{\theta(qt_it_j^{-1};p;q)_{\lambda_i}}
\prod_{j=1}^n\frac{\theta(t_j;p;q)_{|\lambda|}}
{\theta(t_jq^{1+N_j};p;q)_{|\lambda|}}
\\ \nonumber
 && \times \frac{\theta(b,c;p;q)_{|\lambda|}}
{\theta(q/d, q/e;p;q)_{|\lambda|}}
\prod_{j=1}^n \frac{\theta(dt_j, et_j;p;q)_{\lambda_j}}
{\theta(t_jq/b, t_jq/c;p;q)_{\lambda_j}}
\\
&&
=\frac{\theta(q/bd,q/cd;p;q)_{|N|}}{\theta(q/d,q/bcd;p;q)_{|N|}}
\prod_{j=1}^n\frac{\theta(t_jq,t_jq/bc;p;q)_{N_j}}
{\theta(t_jq/b,t_jq/c;p;q)_{N_j}},
\label{e-milne}\end{eqnarray}
where $|\lambda|=\lambda_1+\cdots+\lambda_n$,
$|N|=N_1+\cdots+N_n$, and $bcde=q^{1+|N|}$. For $n=1$ this is the
Frenkel-Turaev sum \cite{fre-tur:elliptic}, and for $p=0$ it is
reduced to the Milne's multiple $_8\Phi_7$ sum for the $A_n$ root
system \cite{mil:multidimensional}.
\end{theorem}
\begin{proof}
We scale $t_i$ for $i=1,\ldots, n$
from the region $|t_i|<1$ to $|t_i|>1$, and keep $|t_{n+1}|, |f_j|<1$
together with the condition $|pq|<|AB|$. During this procedure, some
poles of the integrand $\Delta^I(\mathbf{z};A_n)$ in (\ref{int-an})
go out of the unit disk and, on the contrary, some of them cross
over $\mathbb{T}$ entering inside. The outgoing poles are located
at the following points: $z_k=\{t_iq^{\lambda_i}, i=1,\ldots,n\}$
for each $k=1,\ldots, n$, and the number of such poles is determined
by the conditions $|t_iq^{\lambda_i}|>1$. The ingoing poles
correspond to the points $z_1\cdots z_n=\{t_i^{-1}q^{-\lambda_i},
i=1,\ldots,n\}$.

We denote by $C$ a deformed contour of integration such that none of the
poles mentioned above crosses over $C$ during the change of parameters.
By analyticity, the value of the integral (\ref{int-an}) is not changing
when $C$ replaces $\mathbb{T}$, that is the right-hand side of
(\ref{int-an}) remains the same. If we start to deform the contour $C$ back
to $\mathbb{T}$, we start to pick up residues from the poles by the
Cauchy theorem. As a result, the following formula arises:
\begin{equation}
\int_{C^n}\Delta^I(\mathbf{z};A_n)\prod_{k=1}^n\frac{dz_k}{z_k}
=\sum_{j=0}^{n}\int_{\mathbb{T}^j} R_j(z_1,\ldots,z_j)\prod_{k=1}^j
\frac{dz_k}{z_k},
\label{res}\end{equation}
where $R_n=\Delta^I(\mathbf{z};A_n)$, and $R_j(z_1,\ldots,z_j)$
for $j<n$ are sums of the residues of $\Delta^I(\mathbf{z};A_n)$
corresponding to the poles crossing $C$.

We shall not derive explicit expressions for all coefficients $R_j$
as it was done for the $C_n$ integrals in \cite{die-spi:elliptic,
die-spi:modular}. For our purposes, it suffices to pick up only
the residues that diverge in the limits $f_j\to q^{-N_j}t_j^{-1}$,
$N_j\in \mathbb{N}$, for all $j=1,\ldots, n$ simultaneously.
First, consider the residues appearing from the poles
$z_j= t_jq^{\lambda_j}$, where $\lambda_j$ are some integers such
that $|t_jq^{\lambda_j}|>1$. Straightforward computations yield
\begin{eqnarray}\nonumber
&& R_0^{div}(\mathbf{\lambda})\equiv\prod_{j=1}^n
\lim_{z_j\to t_jq^{\lambda_j}}
(1-t_jq^{\lambda_j}z_j^{-1})\; \Delta^I(\mathbf{z};A_n)
\\ \nonumber && \makebox[2em]{}
= \prod_{\stackrel{i,j=1}{i\neq j}}^n\Gamma^{-1}(t_it_j^{-1}
q^{\lambda_i-\lambda_j})
\frac{\prod_{i=1}^{n+1}\Gamma(t_iDq^{|\lambda|})
\prod_{j=1}^{n+2}\Gamma(f_jD^{-1}q^{-|\lambda|})}
{\prod_{k=1}^n\Gamma(t_kq^{\lambda_k}Dq^{|\lambda|},
D^{-1}q^{-|\lambda|}t_k^{-1}q^{-\lambda_k})}
\\ \nonumber  && \makebox[4em]{} \times
\frac{\prod_{k=1}^n\left(\prod_{\stackrel{i=1}{i\neq k}}^{n+1}
\Gamma(t_it_k^{-1}q^{-\lambda_k}) \prod_{j=1}^{n+2}
\Gamma(f_jt_kq^{\lambda_k})\right)}
{\prod_{k=1}^n\Gamma(AB t_kq^{\lambda_k})\, \Gamma(ABD^{-1}q^{-|\lambda|})}
\\ && \makebox[4em]{} \times
\prod_{k=1}^n\frac{(-1)^{\lambda_k}q^{\lambda_k(\lambda_k+1)/2}}
{(q;q)_\infty(p;p)_\infty\theta(q;p;q)_{\lambda_k}},
\label{residue}\end{eqnarray}
where $D=A/t_{n+1}$. The factors $\Gamma(f_jt_jq^{\lambda_j})$
provide the required divergence in the limits $f_j\to q^{-N_j}t_j^{-1}$.
We write $R_0^{div}(\mathbf{\lambda})= \kappa_n
\Delta(\mathbf{\lambda};A_n)$, where
\begin{eqnarray*}
&& \kappa_n = \prod_{1\leq i<j\leq n}\Gamma^{-1}(t_it_j^{-1},
t_i^{-1}t_j)\frac{\prod_{i=1}^{n+1}\Gamma(t_iD)
\prod_{j=1}^{n+2}\Gamma(f_jD^{-1})}
{\prod_{k=1}^n\Gamma(t_kD,t_k^{-1}D^{-1})}
\\ && \makebox[2em]{} \times
\frac{\prod_{k=1}^n\left(\prod_{\stackrel{i=1}{i\neq k}}^{n+1}
\Gamma(t_it_k^{-1}) \prod_{j=1}^{n+2}\Gamma(f_jt_k)\right)}
{(q;q)_\infty^n(p;p)_\infty^n\prod_{k=1}^n\Gamma(AB t_k)\,\Gamma(ABD^{-1})},
\end{eqnarray*}
and after a chain of simplifying
calculations, $\Delta(\mathbf{\lambda};A_n)$ takes the form
\begin{eqnarray*}
&& \Delta(\mathbf{\lambda};A_n)=
q^{\sum_{j=1}^nj\lambda_j}\prod_{1\leq i<j\leq n}
\frac{\theta(t_it_j^{-1}q^{\lambda_i-\lambda_j};p)}
{\theta(t_it_j^{-1};p)}
\\ && \makebox[4em]{} \times
\prod_{j=1}^n\left(\frac{\theta(t_jD;p;q)_{|\lambda|}}
{\theta(qf_j^{-1}D;p;q)_{|\lambda|}}
\frac{\theta(t_jDq^{|\lambda|+\lambda_j};p)}{\theta(t_jD;p)}
\prod_{k=1}^n\frac{\theta(f_jt_k;p;q)_{\lambda_k}}
{\theta(qt_j^{-1}t_k;p;q)_{\lambda_k}}\right)
\\ && \makebox[4em]{} \times
\frac{\theta(cD,qD/AB;p;q)_{|\lambda|}}
{\theta(qD/d,qD/e;p;q)_{|\lambda|}}
\prod_{k=1}^n\frac{\theta(dt_k,et_k;p;q)_{\lambda_k}}
{\theta(qt_k/c,ABt_k;p;q)_{\lambda_k}},
\end{eqnarray*}
where we have denoted $c=t_{n+1},d=f_{n+1},e=f_{n+2}$.
Now we substitute $f_j=q^{-N_j}t_j^{-1}$ (which assumes that
$AB=q^{-|N|}cde$) in this expression and introduce
the parameter $b\equiv q^{1+|N|}/cde$. As a result,
the function $\Delta(\mathbf{\lambda};A_n)$ becomes equal
to the summand on the left-hand side of (\ref{e-milne}) after
the transformations $t_j\to t_j/D$, $b\to b/D,$ $c\to c/D$,
$d\to Dd$, $e\to De$.

Now, we find the total number of residues of such type.
There is permutational symmetry between the variables $z_1,\ldots, z_n$.
Therefore, there are $n!$ ways to satisfy the equalities
$z_k=t_jq^{\lambda_j}$ using each $t_j$ only once.
The residues of the other outgoing poles located at
$z_k=\{t_iq^{\lambda_i}, i=1,\ldots,n\}$, where at least one $t_i$
enters twice, do not diverge at $f_j\to q^{-N_j}t_j^{-1}$ for some $j$.

We pass to the ingoing poles. It is not difficult to verify that
the residues of $\Delta^I(\mathbf{z};A_n)$ for the poles located
at $z_n=t_j^{-1}q^{-\lambda_j}/z_1\cdots z_{n-1}$ for some
fixed $j$ (or, equivalently, for $z_{n+1}=t_jq^{\lambda_j}$)
are equal to the residues for the poles at $z_n=t_jq^{\lambda_j}$.
Among the remaining poles in the variables $z_1,\ldots,z_{n-1}$,
we must consider only the outgoing ones since only they may
diverge as
$f_k\to q^{-N_k}t_k^{-1}$ with $k=1,\ldots, n, k\neq j$. There are $n$
ways to fix the variable $z_k$ for which we shall consider ingoing
poles, there are $n$ ways to fix the parameter $t_j$ in the equation
$z_{n+1}=t_jq^{\lambda_j}$, and there are $(n-1)!$ appropriate
outgoing poles with the required residue divergence. As a result,
the contribution of these combined ingoing and outgoing poles is equal
to $n^2 (n-1)!$, and the total number of diverging residues
(\ref{residue}) is equal to $(n+1)!$. Roughly speaking, the incoming
poles imitate the $(n+1)$st independent contour of integration over
$z_{n+1}$, which enters symmetrically with $z_1,\ldots,z_n$, and there
are $(n+1)!$ ways to order these variables in the residue calculus.

As has already been mentioned, (\ref{res}) is equal to the right-hand side
of (\ref{int-an}). Now we divide both these expressions by
$(n+1)!\kappa_n$ and take the limits as $f_j\to q^{-N_j}t_j^{-1},
j=1,\ldots,n.$ Since $\kappa_n\to \infty$ in this limit, only the
residues considered above survive in (\ref{res}),
and their sum is given by the elliptic Milne series (\ref{e-milne}).
As to the right-hand side, we get
\begin{eqnarray*}
&& \lim_{f_j\to q^{-N_j}t_j^{-1}}\frac{\text{r.h.s. of (\ref{int-an})}}
{(n+1)!\kappa_n} = \frac{\theta(q/bd,q/cd;p;q)_{|N|}}
{\theta(qD/d,q/Dbcd;p;q)_{|N|}}
\\ && \makebox[4em]{} \times
\prod_{j=1}^n\frac{\theta(qt_jD,qt_j/Dbc;p;q)_{N_j}}
{\theta(qt_j/c,qt_j/b;p;q)_{N_j}},
\end{eqnarray*}
which coincides with the right-hand side of (\ref{e-milne}) after
the appropriate changes of parameters indicated above. The theorem
is proved.
\end{proof}

Formula (\ref{int-an}) generates the following symmetry transformation for
integrals:
\begin{eqnarray}\nonumber
&& \prod_{j=1}^{n+2}\frac{\Gamma(Bf_j^{-1})}{\Gamma(t^{n+1}Bf_j^{-1})}
\int_{\mathbb{T}^n}\frac{\prod_{k=1}^{n+1}\prod_{j=1}^{n+2}
\Gamma(tf_jz_k^{-1},s_jz_k)\prod_{j=1}^n dz_j/z_j}
{\prod_{\stackrel{i,j=1}{i\neq j}}^{n+1}
\Gamma(z_iz_j^{-1})\prod_{k=1}^{n+1}\Gamma(t^{n+1}Sz_k,tBz_k^{-1})}
\\ &&
=\prod_{j=1}^{n+2}\frac{\Gamma(Ss_j^{-1})}{\Gamma(t^{n+1}Ss_j^{-1})}
\int_{\mathbb{T}^n}\frac{\prod_{k=1}^{n+1}\prod_{j=1}^{n+2}
\Gamma(ts_jz_k^{-1},f_jz_k)\prod_{j=1}^n dz_j/z_j}
{\prod_{\stackrel{i,j=1}{i\neq j}}^{n+1}
\Gamma(z_iz_j^{-1})\prod_{k=1}^{n+1}\Gamma(t^{n+1}Bz_k,tSz_k^{-1})}.
\label{an-trans}\end{eqnarray}
Here $t, f_j, s_j, j=1, \ldots, n+2,$ are free independent variables,
$B=\prod_{j=1}^{n+2}f_j, S=\prod_{j=1}^{n+2}s_j$, and it is assumed that
$|t|,|f_j|, |s_j|<1,$ $|pq|<|t^{n+1}B|, |t^{n+1}S|.$
In order to derive this identity, it is necessary to consider the
$2n$-tuple integral
\begin{eqnarray*}
&& \frac{1}{(2\pi i)^n} \int_{\mathbb{T}^{2n}}
\frac{\prod_{k=1}^{n+1}\prod_{j=1}^{n+2}\Gamma(f_jz_k,s_jw_k^{-1})
\prod_{i,k=1}^{n+1}\Gamma(tz_k^{-1}w_i)}
{\prod_{\stackrel{i,j=1}{i\neq j}}^{n+1}
\Gamma(z_iz_j^{-1},w_iw_j^{-1})\prod_{k=1}^{n+1}\Gamma(t^{n+1}Bz_k,
t^{n+1}Sw_k^{-1})}
\\ && \makebox[4em]{} \times
\frac{dz_1}{z_1}\cdots\frac{dz_n}{z_n}
\frac{dw_1}{w_1}\cdots\frac{dw_n}{w_n},
\end{eqnarray*}
where $z_1\cdots z_{n+1}=w_1\cdots w_{n+1}=1$.
Integration with respect to the variables $z_k$ with the help of
(\ref{int-an}) makes this expression proportional to the left-hand
side of (\ref{an-trans}) (after the replacements $w_k\to z_k^{-1}$).
Changing the order of integration, which is allowed because the integrand
is bounded on $\mathbb{T}$, we arrive at the integral standing on the
right-hand side of (\ref{an-trans}), up to some coefficient.
After cancelling common factors, we get the required identity.
For $p=0$ this reduces to the Denis-Gustafson
transformation formula \cite{den-gus:beta}, which describes a Bailey
transformation for a terminating multivariable $_{10}\Phi_9$ series.
It is natural to expect that in our case formula (\ref{an-trans})
yields a Bailey transformation for some terminating $A_n$ multiple
$_{12}E_{11}$ series appearing from sums of residues of the
corresponding integrals.

\section{Some other $A_n$ integrals}

Now from conjecture (\ref{int-an}) we derive a number of
different multiple $A_n$-integrals. Denote
\begin{eqnarray}\nonumber
&& \Delta^{II}({\bf z};A_n)=
\frac{1}{(2\pi i)^n}\prod_{1\leq i<j\leq n+1}
\frac{\Gamma(tz_iz_j,sz_i^{-1}z_j^{-1})}
{\Gamma(z_iz_j^{-1},z_i^{-1}z_j)}
\\ && \makebox[4em]{} \times
\prod_{j=1}^{n+1}\frac{\Gamma(t_1z_j,t_2z_j,t_3z_j,t_4z_j^{-1},t_5z_j^{-1})}
{\Gamma(z_j(ts)^{n-1}\prod_{j=1}^5 t_j)}.
\label{a2}\end{eqnarray}

\begin{theorem}
Suppose the validity of the conjectured $A_n$ and $C_n$ multiple
elliptic beta integrals (\ref{int-an}) and (\ref{SintA}), respectively.
Then, the following two integration formulae are true.
For odd $n=2m-1$, we have
\begin{eqnarray}\nonumber
&& \int_{\mathbb{T}^n}\Delta^{II}({\bf z};A_n)
\frac{dz_1}{z_1}\ldots \frac{dz_n}{z_n}
= \frac{ (n+1)! }{ (q;q)_\infty^n (p;p)_\infty^n }
\\ && \makebox[4em]{} \times
\frac{ \Gamma(t^m,s^m,s^{m-1}t_4t_5)\prod_{1\leq i<j\leq 3}
\Gamma(t^{m-1}t_it_j) }
{ \prod_{k=4}^5\Gamma(t^{2m-2}s^{m-1}t_1t_2t_3t_k) }
\nonumber \\ && \makebox[4em]{} \times
\prod_{j=1}^m \frac{ \prod_{i=1}^3\prod_{k=4}^5
\Gamma((ts)^{j-1}t_it_k) }
{ \prod_{1\leq i<\ell\leq 3} \Gamma((ts)^{m+j-2}t_it_\ell t_4t_5) }
\nonumber \\ && \makebox[4em]{} \times
\prod_{j=1}^{m-1} \frac{ \Gamma((ts)^j,t^js^{j-1}t_4t_5)
\prod_{1\leq i<\ell\leq 3}\Gamma(t^{j-1}s^jt_it_\ell) }
{ \prod_{k=4}^5\Gamma(t^{m+j-2}s^{m+j-1}t_1t_2t_3 t_k) }.
\label{int-an-odd}\end{eqnarray}

For even $n=2m$, we have
\begin{eqnarray}\nonumber
&&\int_{\mathbb{T}^n}\Delta^{II}({\bf z};A_n)
\frac{dz_1}{z_1}\ldots \frac{dz_n}{z_n}
= \frac{ (n+1)! }{ (q;q)_\infty^n (p;p)_\infty^n }
\\ && \makebox[4em]{} \times
\frac{ \prod_{i=1}^3\Gamma(t^mt_i)\prod_{k=4}^5\Gamma(s^mt_k)\,
\Gamma(t^{m-1}t_1t_2t_3) }
{ \Gamma(t^{2m-1}s^{m-1}\prod_{i=1}^5t_i,t^{2m-1}s^m t_1t_2t_3) }
\nonumber \\ && \makebox[4em]{} \times
\prod_{j=1}^m \frac{ \Gamma((ts)^j,t^js^{j-1}t_4t_5)
\prod_{i=1}^3\prod_{k=4}^5 \Gamma((ts)^{j-1}t_it_k) }
{ \prod_{k=4}^5\Gamma(t^{m+j-2}s^{m+j-1}t_k^{-1}\prod_{i=1}^5t_i)}
\nonumber \\ && \makebox[4em]{} \times
\prod_{j=1}^m \prod_{1\leq i<\ell\leq3 }
\frac{\Gamma(t^{j-1}s^jt_it_\ell)}
{\Gamma((ts)^{m+j-1}t_it_\ell t_4t_5)}.
\label{int-an-even}\end{eqnarray}
\end{theorem}

\begin{proof}
The proofs follow the procedure used by Gustafson in \cite{gus:some2}
for proving the $p=0$ cases of the integrals (\ref{int-an-odd}) and
(\ref{int-an-even}). We start with the case of odd $n=2m-1$.
Consider the following $(4m-1)$-tuple integral:
\begin{eqnarray}\nonumber
\lefteqn{ \int_{\mathbb{T}^{4m-1}} \frac{\prod_{i=1}^{2m}\prod_{j=1}^m
\Gamma(t^{1/2}z_iw_j,t^{1/2}z_iw_j^{-1},
s^{1/2}z_i^{-1}x_j,s^{1/2}z_i^{-1}x_j^{-1})}
{\prod_{i,j=1;\: i\neq j}^{2m}\Gamma(z_iz_j^{-1}) \prod_{\nu=\pm1}
\prod_{1\leq i<j\leq m}\Gamma(w_i^\nu w_j^\nu,w_i^\nu w_j^{-\nu},
x_i^\nu x_j^\nu,x_i^\nu x_j^{-\nu})} }&&
\\ && \times
\prod_{i=1}^{2m}\frac{\Gamma(z_i(ts)^{m-2}\prod_{k=1}^5 t_k)}
{\Gamma(z_i(ts)^{2m-2}\prod_{k=1}^5 t_k)}
\prod_{k=1}^{2m-1}\frac{dz_k}{z_k}\prod_{j=1}^m\Biggl(\frac{dw_j}{w_j}
\frac{dx_j}{x_j} \label{n=2m-1} \\ && \times
\prod_{\nu=\pm1} \Biggl( \frac{\prod_{k=1}^3\Gamma(t^{-1/2}t_kw_j^\nu)
\Gamma(x_j^\nu t^{m-2}s^{-1/2}t_1t_2t_3)
\prod_{k=4}^5\Gamma(s^{-1/2}t_kx_j^\nu)}
{\Gamma(w_j^\nu t^{m-3/2}t_1t_2t_3,x_j^\nu t^{m-2}s^{m-3/2}\prod_{k=1}^5t_k,
w_j^{2\nu},x_j^{2\nu})} \Biggr)\Biggr),
\nonumber\end{eqnarray}
where $\prod_{i=1}^{2m}z_j=1$.
Using the exact $C_n$ integration formula of type $I$ (see (\ref{SintA})),
first we take integrals in (\ref{n=2m-1}) with respect
to the variables $w_j, j=1,\ldots,m,$ and after that with respect to
$x_j, j=1,\ldots,m$. The resulting integral is equal to the left-hand
side of (\ref{int-an-odd}) up to the factor
\begin{eqnarray*}
&& (2\pi i)^{4m-1}\frac{2^{2m}(m!)^2}{(p;p)_\infty^{2m}(q;q)_\infty^{2m}}
\frac{\Gamma(s^{-1}t_4t_5)}{\Gamma(s^{m-1}t_4t_5)}
\\ && \times
\prod_{1\leq i<k\leq 3}\frac{\Gamma(t^{-1}t_it_k)}{\Gamma(t^{m-1}t_it_k)}
\prod_{k=4,5}\frac{\Gamma(t^{m-2}s^{-1}t_k^{-1}\prod_{i=1}^5t_k)}
{\Gamma(t^{m-2}s^{m-1}t_k^{-1}\prod_{i=1}^5t_k)}.
\end{eqnarray*}

In this two step procedure, we need the following restrictions upon the parameters:
$$
|t|<1, \quad |t_{1,2,3}|<|t|^{1/2}, \quad |pq|<|t^{m-3/2}t_1t_2t_3|
$$
and
$$
|s|<1,\quad |t_{4,5}|<|s|^{1/2},\quad
|pq|<|t^{m-2}s^{m-3/2}\prod_{k=1}^5t_k|,
$$
respectively. However, the resulting expression can be extended
analytically to the region $|t_k|<1, \, k=1,\ldots,5$, $|pq|<|(ts)^{2m-2}
\prod_{k=1}^5t_k|$ without changing the integral value.

Since the integrand function in (\ref{n=2m-1}) is bounded on the unit
circle, we can change the order of integrations. First, we take the
integrals over $z_i,\: i=1,\ldots,2m-1,$ using the $A_n$-formula
(\ref{int-an}). Then we apply formula (\ref{SintA}) in order to take
the integrals over $x_j, j=1,\ldots,m$. Finally, we apply the intrinsic
elliptic Selberg integral (\ref{SintB}) for taking the integrals over
$w_j,\: j=1,\ldots,m;$ this leads to the following expression:
\begin{eqnarray*}
&& \frac{(2\pi i)^{4m-1}(2m)! (2^mm!)^2}{((p;p)_\infty(q;q)_\infty)^{4m-1}}
\frac{\Gamma(t^m,s^m,t_4t_5s^{-1})} {\Gamma((ts)^m,t^ms^{m-1}t_4t_5)}
\prod_{k=4}^5\frac{\Gamma(t^{m-2}s^{-1}t_k^{-1}\prod_{i=1}^5t_i)}
{\Gamma(t^{2m-2}s^{m-1}t_k^{-1}\prod_{i=1}^5t_i)}
\\ && \times
\prod_{j=1}^m\frac{\Gamma((ts)^j,t^js^{j-1}t_4t_5)
\prod_{1\leq i<k\leq 3}\Gamma(t^{j-2}s^{j-1}t_it_k)
\prod_{k=4}^5 \prod_{i=1}^3\Gamma((ts)^{j-1}t_kt_i)}
{\prod_{k=1}^3\Gamma((ts)^{m+j-2}t_k^{-1}\prod_{i=1}^5t_i)
\prod_{k=4}^5 \Gamma(t^{m+j-3}s^{m+j-2}t_k^{-1}\prod_{i=1}^5t_i)}.
\end{eqnarray*}
As a result, we get the needed integration formula for odd $n=2m-1$.

In order to prove (\ref{int-an-even}), we consider the $4m$-tuple integral
\begin{eqnarray}\nonumber
&& \int_{\mathbb{T}^{4m}} \frac{\prod_{i=1}^{2m+1}\prod_{j=1}^m
\Gamma(t^{1/2}z_iw_j,t^{1/2}z_iw_j^{-1},
s^{1/2}z_i^{-1}x_j,s^{1/2}z_i^{-1}x_j^{-1})}
{\prod_{i,j=1;\: i\neq j}^{2m+1}\Gamma(z_iz_j^{-1}) \prod_{\nu=\pm1}
\prod_{1\leq i<j\leq m}\Gamma(w_i^\nu w_j^\nu,w_i^\nu w_j^{-\nu},
x_i^\nu x_j^\nu,x_i^\nu x_j^{-\nu})}
\\ && \makebox[2em]{}\times
\prod_{i=1}^{2m+1}\frac{\Gamma(t_3z_i,s^{m-1}t_4t_5z_i,
t^{m-1}t_1t_2z_i^{-1})}{\Gamma(z_i(ts)^{2m-1}\prod_{k=1}^5 t_k)}
\prod_{j=1}^m\Biggl( \frac{dw_j}{w_j}\frac{dx_j}{x_j}
\label{n=2m} \\ && \makebox[2em]{}\times
\prod_{\nu=\pm1} \Biggl(\frac{\prod_{k=1}^2\Gamma(t^{-1/2}t_kw_j^\nu,
s^{-1/2}t_{k+3}x_j^\nu)}{\Gamma(t^{m-1/2}t_1t_2w_j^\nu,
s^{m-1/2}t_4t_5x_j^\nu, w_j^{2\nu},x_j^{2\nu})} \Biggr)\Biggr)
\prod_{k=1}^{2m}\frac{dz_k}{z_k},
\nonumber\end{eqnarray}
where $\prod_{i=1}^{2m+1}z_j=1$. Repeating the same trick as in the
case of odd $n$ (that is, integrating successively with respect to
the variables $w_j$ and $x_j$ and then changing the order of
integrations in this expression), we get (\ref{int-an-even}).
\end{proof}

In a similar way, we can establish elliptic analogs of the $A_n$
basic hypergeometric integrals of Gustafson and Rakha \cite{gus-rak:beta}.

\begin{theorem}
Suppose the validity of the $A_n$ and $C_n$ multiple elliptic beta
integrals (\ref{int-an}) and (\ref{SintA}), respectively. Denote
\begin{eqnarray}\nonumber
&& \Delta^{III}({\bf z};A_n)=
\frac{1}{(2\pi i)^n}\prod_{1\leq i<j\leq n+1}
\frac{\Gamma(tz_iz_j)}{\Gamma(z_iz_j^{-1},z_i^{-1}z_j) }
\\ && \makebox[4em]{} \times
\frac{\prod_{i=1}^{n+1}\left(\prod_{k=1}^{n+1}\Gamma(t_kz_i^{-1})
\prod_{k=n+2}^{n+4}\Gamma(tt_kz_i)\right)}
{\prod_{j=1}^{n+1}\Gamma(Az_j^{-1}) },
\label{a3}\end{eqnarray}
where $A=t^{n+2}\prod_{i=1}^{n+4}t_i$ and $\prod_{j=1}^{n+1}z_j=1$.
Then the following two integration formulae are true.
For odd $n=2l-1$, we have
\begin{eqnarray}
&& \int_{\mathbb{T}^n}\Delta^{III}({\bf z};A_n)
\frac{dz_1}{z_1}\ldots\frac{dz_n}{z_n}
= \frac{ (n+1)! }{ (q;q)_\infty^n (p;p)_\infty^n }
\frac{\Gamma(t^l,\prod_{k=1}^{2l}t_k)}
{\Gamma(t^l\prod_{k=1}^{2l}t_k)}
\label{int-an-o} \\ && \makebox[1em]{} \times
\frac{ \prod_{i=1}^{2l}\prod_{j=2l+1}^{2l+3}\Gamma(tt_it_j)
\prod_{1\leq i<j\leq 2l} \Gamma(tt_it_j)
\prod_{2l+1\leq i<j\leq 2l+3}\Gamma(t^{l+1}t_it_j)}
{\prod_{i=1}^{2l}\Gamma(t^{2l+1}t_i^{-1}\prod_{k=1}^{2l+3}t_k)
\prod_{i=2l+1}^{2l+3}\Gamma(t^{l+1}t_i^{-1}\prod_{k=1}^{2l+3}t_k) }.
\nonumber\end{eqnarray}

For even $n=2l$, we have
\begin{eqnarray}\nonumber
&&\int_{\mathbb{T}^n}\Delta^{III}({\bf z};A_n)
\frac{dz_1}{z_1}\ldots\frac{dz_n}{z_n}
= \frac{ (n+1)! }{ (q;q)_\infty^n (p;p)_\infty^n }
\frac{\Gamma(\prod_{k=1}^{2l+1}t_k,t^{l+2}\prod_{k=2l+2}^{2l+4}t_k)}
{\Gamma(t^{l+2}\prod_{k=1}^{2l+4}t_k)}
\\ && \makebox[1em]{} \times
\frac{ \prod_{i=1}^{2l+1}\prod_{j=2l+2}^{2l+4}\Gamma(tt_it_j)
\prod_{1\leq i<j\leq 2l+1} \Gamma(tt_it_j)
\prod_{i=2l+2}^{2l+4}\Gamma(t^{l+1}t_i)}
{\prod_{i=1}^{2l+1}\Gamma(t^{2l+2}t_i^{-1}\prod_{k=1}^{2l+4}t_k)
\prod_{i=2l+2}^{2l+4}\Gamma(t^{l+1}t_i\prod_{k=1}^{2l+1}t_k) }.
\label{int-an-e} \end{eqnarray}
\end{theorem}

\begin{proof}
In accordance with the procedure used in \cite{gus-rak:beta},
we consider the $(3l-1)$-tuple integral
\begin{eqnarray*}
&& \int_{\mathbb{T}^{3l-1}}
\frac{\prod_{i=1}^{2l}\prod_{j=1}^l\Gamma(t^{1/2}z_iw_j,t^{1/2}z_iw_j^{-1})
\prod_{i=0}^{2l}\prod_{j=1}^{2l}\Gamma(t_iz_j^{-1})}
{\prod_{i,j=1; i\neq j}^{2l}\Gamma(z_iz_j^{-1})
\prod_{j=1}^{2l} \Gamma(t^l\prod_{i=0}^{2l}t_iz_j^{-1})}
\\ && \makebox[1em]{} \times
\prod_{\nu=\pm1}\prod_{1\leq i<j\leq l}\Gamma^{-1}
(w_i^\nu w_j^{\nu},w_i^\nu w_j^{-\nu})
\prod_{j=1}^l
\frac{\prod_{k=2l+1}^{2l+3}\Gamma(t^{1/2}t_kw_j^\nu)}
{\Gamma(w_j^{2\nu},t^{l+3/2}\prod_{k=2l+1}^{2l+3}t_kw_j^\nu)}
\\ && \makebox[1em]{} \times
\frac{dw_1}{w_1}\ldots\frac{dw_l}{w_l}\frac{dz_1}{z_1}
\ldots\frac{dz_{2l-1}}{z_{2l-1}},
\end{eqnarray*}
where $\prod_{i=1}^{2l}z_i=1$ and $t_0=t^{l+1}\prod_{k=2l+1}^{2l+3}t_k$.
Integrating with respect to the variables $w_j$ with the help of
formula (\ref{SintA}), we get the left-hand side of (\ref{int-an-o})
up to some factor. Changing the order of
integration, we can integrate over $z_i$ using (\ref{int-an})
(where it is necessary to change $z_k$ to $z_k^{-1}$)
and then over $w_j$ using (\ref{SintA}). Equating two expressions,
we arrive at formula (\ref{int-an-o}).

In a similar way, in the case of even $n=2l$ we consider the $(3l+1)$-tuple
integral
\begin{eqnarray*}
&& \int_{\mathbb{T}^{3l+1}}
\frac{\prod_{i=1}^{2l+1}\left(\prod_{j=1}^{l+1}
\Gamma(t^{1/2}z_iw_j,t^{1/2}z_iw_j^{-1})
\prod_{j=1}^{2l+1}\Gamma(t_iz_j^{-1})\right)}
{\prod_{i,j=1; i\neq j}^{2l+1}\Gamma(z_iz_j^{-1})
\prod_{j=1}^{2l+1} \Gamma(t^{l+1}\prod_{i=1}^{2l+1}t_iz_j)}
\\ && \makebox[1em]{} \times
\prod_{\nu=\pm1}\prod_{1\leq i<j\leq l+1}\Gamma^{-1}
(w_i^\nu w_j^{\nu},w_i^\nu w_j^{-\nu})
\prod_{j=1}^{l+1}
\frac{\prod_{k=2l+1}^{2l+5}\Gamma(t^{1/2}t_kw_j^\nu)}
{\Gamma(w_j^{2\nu},t^{l+5/2}\prod_{k=2l+1}^{2l+5}t_kw_j^\nu)}
\\ && \makebox[1em]{} \times
\frac{dw_1}{w_1}\ldots\frac{dw_{l+1}}{w_{l+1}}\frac{dz_1}{z_1}
\ldots\frac{dz_{2l}}{z_{2l}},
\end{eqnarray*}
where $\prod_{i=1}^{2l+1}z_i=1$ and $t_{2l+5}=t^l\prod_{k=1}^{2l+1}t_k$.
Repeating the same trick as in the preceding case, we get the desired
integration formula (\ref{int-an-e}).
\end{proof}

Sums of residues for the derived integrals
(\ref{int-an-odd})-(\ref{int-an-e}) form elliptic hypergeometric
series on the $A_n$ root system that differ from the series (\ref{e-milne})
introduced in \cite{ros:elliptic,spi:modularity}. We skip their
consideration and formulate only a conjecture concerning the
elliptic extension of Theorem 1.2 in \cite{gus-rak:beta}.

\begin{conjecture}
Suppose that $N$ is a positive integer and $\prod_{k=1}^n t_k=q^{-N}$.
Then
\begin{eqnarray}\nonumber
&&
\sum_{\stackrel{\lambda_k=0,\ldots,N}{\lambda_1+\ldots+\lambda_n=N}}
\frac{
\prod_{1\leq i<j\leq n}\theta(tt_it_j)_{\lambda_i+\lambda_j}
\prod_{i=1}^n\prod_{j=n+1}^{n+3}\theta(tt_it_j)_{\lambda_i}
\prod_{i,j=1}^n\theta(t_it_j^{-1})_{-\lambda_j}
}{
\prod_{i,j=1;i\neq j}^n\theta(t_it_j^{-1})_{\lambda_i-\lambda_j}
\prod_{j=1}^n\theta(t^{n+1}t_j^{-1}\prod_{k=1}^{n+3}t_k)_{-\lambda_j}
} \\ && \makebox[2em]{}
= \left\{ \begin{aligned} % \begin{array}{ll}
\frac{\theta(1)_{-N}}{\theta(t^{n/2})_{-N}
\prod_{n+1\leq i<j\leq n+3}\theta(t^{(n+2)/2}t_it_j)_{-N}},
& \quad n \quad \text{is even,}  \\
\frac{\theta(1)_{-N}}{\prod_{i=n+1}^{n+3}\theta(t^{(n+1)/2}t_i)_{-N}
\theta(t^{(n+3)/2}\prod_{i=n+1}^{n+3}t_i)_{-N} },
& \quad n \quad \text{is odd,}  \\
\end{aligned} \right.
\label{new-sums} \end{eqnarray}
where $\theta(a)_\lambda\equiv\theta(a;p;q)_\lambda$.
\end{conjecture}

These summation formulas are expected to follow from the residue calculus
for the integrals (\ref{int-an-o}) and (\ref{int-an-e}).
Some evidence in favor of conjecture (\ref{new-sums})
is provided by the following theorem.

\begin{theorem}
Denote $t=q^g, t_i=q^{g_i},i=1,\ldots,n+3$ (so that
$\sum_{j=1}^{n}g_j+N=0$). The series
$\sum_{\mathbf{\lambda}}c(\mathbf{\lambda})$
standing on the left-hand side of (\ref{new-sums}) is a totally
elliptic hypergeometric series, that is, the ratios of successive
series coefficients
$$
h_k(\mathbf{\lambda})=\frac{c(\lambda_1,\ldots,\lambda_k+1,\ldots,
\lambda_n)}{c(\lambda_1,\ldots,\lambda_n)}
$$
are elliptic functions of all unconstrained variables in the set
$(\lambda_1,\ldots,\lambda_n, g,g_1,\ldots,$ $ g_{n+3})$.
Moreover, the functions $h_k(\mathbf{\lambda})$ are $SL(2,\mathbb{Z})$
modular invariant. The ratios of the expressions standing on the two sides
of (\ref{new-sums}) are elliptic functions of $g$ and $n+2$ free
parameters in the set $(g_1,\ldots,g_{n+3})$, and these ratios
are modular invariant as well.
\end{theorem}

We skip the proof of this theorem, because it consists of quite long but
straightforward computations whose structure was described in detail
in \cite{spi:theta,spi:modularity} in the process of similar considerations
for different elliptic hypergeometric series summation formulas. Using
the fact that there are no cusp forms of weights below 12,
as in \cite{die-spi:elliptic,spi:modularity} from this theorem
we deduce that relations (\ref{new-sums}) are valid in the small
$\sigma$ expansion up to the terms of order of $\sigma^{12}$.
This gives also yet another example in favor of the general conjecture
of \cite{spi:theta} that all totally elliptic hypergeometric series are
automatically modular invariant.

The integrals (\ref{int-an-odd}) and (\ref{int-an-even}) are
expected to generate $A_n$ summation formulas similar to
(\ref{new-sums}). It is also natural to expect that all multiple elliptic
beta integrals described above lead to integral representations for
various multiple $_{12}E_{11}$ elliptic hypergeometric series
generalizing the single variable formula announced in \cite{spi:special}
(see Appendix B for the proof of it). We suppose that, sas in the
$A_n$-case, there exist several types of elliptic beta
integrals and elliptic hypergeometric series
sums associated with the $D_n$ root system (see, e.g.,
\cite{ros:elliptic,spi:modularity}), but their consideration
lies beyond the scope of the present paper.

\section{Relations to the generalized eigenvalue problems}

Consider the very-well-poised theta hypergeometric series (\ref{vwp-1})
with the additional constraint $(-q)^ne^{P_3(n)}=(qx)^n$ for some
$x\in\mathbb{C}$. Special notation for such series was introduced
in \cite{spi:bailey}:
\begin{eqnarray} \nonumber
\lefteqn{ _{r+1}V_r(t_0;t_1,\ldots,t_{r-4};q,p;x) }&&
\\ && \makebox[3em]{}
\equiv \sum_{n=0}^\infty \frac{\theta(t_0q^{2n};p)}{\theta(t_0;p)}
\prod_{m=0}^{r-4}\frac{\theta(t_m;p;q)_n}{\theta(qt_0t_m^{-1};p;q)_n}
\, (qx)^n.
\label{vwp-2}\end{eqnarray}
For (\ref{vwp-2}), the balancing condition (\ref{theta-balance}) is
reduced to
$$
\prod_{m=1}^{r-4}t_m= t_0^{(r-5)/2}q^{(r-7)/2}.
$$
As $p\to 0$, the $_{r+1}V_r$ series are reduced to $_{r-1}W_{r-2}$
very-well-poised $q$-hypergeome\-tric series of the argument $qx$
(in the notation of \cite{gas-rah:basic}).

The balanced $_{12}V_{11}$ series with $x=1$ plays an important
role in applications. For instance, the elliptic solutions of the
Yang-Baxter equation derived by Date et al in
\cite{djkmo:exactly1,djkmo:exactly2} are expressed in terms of such
series for a particular choice of parameters \cite{fre-tur:elliptic}.
In all cases of the $_{12}V_{11}$ function to be considered below, we
have $x=1$; therefore, we omit dependence on this unit argument
from now on.

We denote $\mathcal{E}(\mathbf{t}) \equiv
{_{12}}V_{11}(t_0;t_1,\ldots,t_7;q,p),$
where $\prod_{m=1}^7t_m=t_0^3q^2$, and assume that this series terminates
due to the condition $t_m=q^{-n},\, n\in\mathbb{N},$ for some $m$.
In \cite{spi-zhe:spectral,spi-zhe:classical}, the following two
contiguous relations for $\mathcal{E}(\mathbf{t})$ were derived:
\begin{eqnarray}\label{1_con}
\lefteqn{\mathcal{E}(\mathbf{t}) - \mathcal{E}(t_0;t_1,\ldots,t_5,q^{-1}t_6,
qt_7) } &&  \\ &&
=\frac{\theta(qt_0,q^2t_0,qt_7/t_6,t_6t_7/qt_0;p)}
{\theta(qt_0/t_6,q^2t_0/t_6,t_0/t_7,t_7/qt_0;p)}
\prod_{r=1}^5\frac{\theta(t_r;p)}{\theta(qt_0/t_r;p)}\,
\mathcal{E}(q^2t_0;qt_1,\ldots,qt_5,t_6,qt_7),
 \nonumber\\  \nonumber
\lefteqn{\frac{\theta(t_7;p)}{\theta(t_6/qt_0,t_6/q^2t_0,t_6/t_7;p)}
\prod_{r=1}^5 \theta(t_rt_6/qt_0;p)\,
\mathcal{E}(q^2t_0;qt_1,\ldots,qt_5,t_6,qt_7) } &&
\\ && \makebox[4em]{}
+\frac{\theta(t_6;p)}{\theta(t_7/qt_0,t_7/q^2t_0,t_7/t_6;p)}
\prod_{r=1}^5 \theta(t_rt_7/qt_0;p)\,
\mathcal{E}(q^2t_0;qt_1,\ldots,qt_6,t_7)
\nonumber \\ && \makebox[8em]{}
=\frac{1}{\theta(qt_0,q^2t_0;p)}
\prod_{r=1}^5\theta(qt_0/t_r;p)\, \mathcal{E}(\mathbf{t}).
\label{2_con} \end{eqnarray}
Their combination yields
\begin{eqnarray}\nonumber
\lefteqn{
\frac{\theta(t_7,t_0/t_7,qt_0/t_7;p)}{\theta(qt_7/t_6,t_7/t_6;p)}
\prod_{r=1}^5\theta(qt_0/t_6t_r;p)\left(
\mathcal{E}(t_0;t_1,\ldots,t_5,q^{-1}t_6,qt_7) -
\mathcal{E}(\mathbf{t}) \right) } &&
\\  \nonumber
&& +\frac{\theta(t_6,t_0/t_6,qt_0/t_6;p)}
{\theta(qt_6/t_7,t_6/t_7;p)}\prod_{r=1}^5\theta(qt_0/t_7t_r;p)
\left(\mathcal{E}(t_0;t_1,\ldots,t_5,qt_6,q^{-1}t_7)-
\mathcal{E}(\mathbf{t})\right) \\
&& \makebox[4em]{}
+\theta(qt_0/t_6t_7;p)\prod_{r=1}^5\theta(t_r;p)\,
\mathcal{E}(\mathbf{t})=0.
\label{3_con}\end{eqnarray}
As $p\to 0$, these three equalities are reduced to the contiguous
relations for the terminating very-well-poised balanced $_{10}\Phi_9$
series of Gupta and Masson \cite{gup-mas:contiguous}. Similar
contiguous relations at the level of $_8\Phi_7$ functions
were constructed earlier by Ismail and Rahman \cite{ism-rah:associated}.

We change the parametrization of the $\mathcal{E}$-function
and consider relation (\ref{3_con}) for the following function:
\begin{equation}
R_{n}(z;q,p)\equiv {_{12}V_{11}}\left(\frac{t_3}{t_4};
\frac{q}{t_0t_4},\frac{q}{t_1t_4},\frac{q}{t_2t_4},t_3z,\frac{t_3}{z},q^{-n},
\frac{Aq^{n-1}}{t_4};q,p\right),
\label{R_n}\end{equation}
where $A=\prod_{r=0}^4 t_r$. We replace the parameter $t_0$
in (\ref{3_con}) by $t_3/t_4$; the variables $t_1,t_2,$ and $t_3$ are
replaced by $q/t_0t_4,q/t_1t_4,$ and $q/t_2t_4$; the variables
$t_4$ and $t_5$ by $q^{-n}$ and $Aq^{n-1}/t_4$; and
the variables $t_6$ and $t_7$ by $t_3z$ and $t_3/z$,
respectively. As a result, we see that $R_{n}(z;q,p)$ provides a
particular solution of the following finite-difference equation:
\begin{equation}
\mathcal{D}_\mu f(z)=0,\qquad
\mathcal{D}_\mu=V_\mu(z)(T-1)+V_\mu(z^{-1})(T^{-1}-1)+\kappa_\mu,
\label{eig}\end{equation}
where $T$ is the $q$-shift operator, $Tf(z)=f(qz)$, and
\begin{eqnarray}\label{V}
&& V_\mu(z)=\theta\left(\frac{t_4}{q\mu z},\frac{A\mu}{q^2z},\frac{t_4z}{q};p\right)
\frac{\prod_{r=0}^4\theta(t_rz;p)}{\theta(z^2,qz^2;p)},
\\
&& \kappa_\mu=\theta\left(\frac{A\mu}{qt_4},\mu^{-1};p\right)
\prod_{r=0}^3\theta\left(\frac{t_rt_4}{q};p\right).
\label{kappa}\end{eqnarray}
The functions $f(z)=R_n(z;q,p)$ solve (\ref{eig})
for $\mu=q^{n}, n\in\mathbb{N}.$

Equation (\ref{eig}) looks like a nonstandard eigenvalue problem with the
``spectral parameter" $\mu$; indeed, it can be rewritten as the
generalized eigenvalue problem
\begin{equation}\label{gevp}
 \mathcal{D}_\eta f(z)=\lambda \mathcal{D}_\xi f(z),
\end{equation}
where the spectral parameter $\lambda$ is
\begin{equation}
\lambda=\frac{\theta\left(\frac{\mu A\eta}{qt_4},
\frac{\mu }{\eta};p\right)}{\theta\left(\frac{\mu A\xi}{qt_4},
\frac{\mu }{\xi};p\right)}
\label{lambda}\end{equation}
and the operators $\mathcal{D}_{\xi}$, $\mathcal{D}_{\eta}$ are obtained
from $\mathcal{D}_{\mu}$ after the replacements of $\mu$ by
arbitrary gauge parameters $\xi,\eta \in\mathbb{C},\, \xi\neq\eta p^k,
qt_4p^k/A\eta,$ $k\in\mathbb{Z}$.
Application of the theta function identity
(\ref{ident}) to equation (\ref{gevp}) yields
$$
\mathcal{D}_\eta-\lambda \mathcal{D}_\xi
=\frac{\theta(A\eta\xi/qt_4,\xi/\eta;p)}
{\theta(A\mu \xi/qt_4, \xi/\mu ;p)}\, \mathcal{D}_{\mu},
$$
which shows that the gauge parameters $\xi,\eta$ drop out completely
from the equation determining $f(z)$, $\mathcal{D}_\mu f(z)=0$.

A three term recurrence relation for the functions $R_n(z;q,p)$ was
derived in \cite{spi-zhe:spectral,spi-zhe:classical}. It appears
from formula (\ref{3_con}) if there we replace $t_6$ by $q^{-n}$ and
$t_7$ by $Aq^{n-1}/t_4$, and substitute  $t_1\to q/t_0t_4,
t_2\to q/t_1t_4, t_3\to q/t_2t_4$, $t_4\to t_3z$,
$t_5\to t_3/z$. After some work, this foirmula can
be represented in the form
\begin{eqnarray}\nonumber
&& (\gamma(z)-\alpha_{n+1})B(Aq^{n-1}/t_4)\left(R_{n+1}(z;q,p)-R_n(z;q,p)
\right) \\ \nonumber && \makebox[2em]{}
+(\gamma(z)-\beta_{n-1})B(q^{-n})
\left(R_{n-1}(z;q,p)-R_n(z;q,p)\right)
\\ && \makebox[4em]{}
+\delta\left(\gamma(z)-\gamma(t_3)\right)R_n(z;q,p)=0,
\label{ttr}\end{eqnarray}
where
\begin{eqnarray}\label{B}
&& B(x)=\frac{\theta\left(x,\frac{t_3}{t_4x},
\frac{qt_3}{t_4x},\frac{qx}{t_0t_1},\frac{qx}{t_0t_2},
\frac{qx}{t_1t_2},\frac{q^2\eta x}{A},\frac{q^2x}{A\eta};p
\right)}{\theta\left(\frac{qt_4x^2}{A},\frac{q^2t_4x^2}{A};p\right)},
\\ \label{gamma}
&& \delta=\theta\left(\frac{q^2t_3}{A},\frac{q}{t_0t_4},
\frac{q}{t_1t_4},\frac{q}{t_2t_4},t_3\eta,\frac{t_3}{\eta};p\right),
\\ &&
\gamma(z)=\frac{\theta(z\xi,z/\xi;p)}{\theta(z\eta,z/\eta;p)},
\label{ab} \\ &&
\alpha_n=\gamma(q^n/t_4),\qquad \beta_n=\gamma(q^{n-1}A).
\end{eqnarray}
Here $\xi$ and $\eta \neq \xi p^k,\xi^{-1}p^k,$ $k\in\mathbb{Z},$
are arbitrary gauge parameters (they are not related to $\xi,\eta$
in the difference equation, but we use the same notation).
Substituting (\ref{B})-(\ref{ab}) in (\ref{ttr}) and applying
identity (\ref{ident}), we see that the auxiliary gauge parameters
$\xi, \eta$ drop out completely from the resulting recurrence relation.

Since $B(q^{-n})=0$ for $n=0$, the indeterminate $R_{-1}$
is not involved in (\ref{ttr}) for $n=0$. We can say that $R_n(z;q,p)$
are generated by the three term recurrence relation (\ref{ttr})
for the initial conditions $R_{-1}=0,R_0=1$. All recurrence
coefficients in (\ref{ttr}) depend linearly on the
variable $\gamma(z)$, which absorbs $z$-dependence. Therefore,
$R_n(z;q,p)$ are rational functions of $\gamma(z)$
with $n$ being the degree of polynomials in $\gamma(z)$ in
the numerator and denominator of $R_n$. Moreover, the poles of these
functions are located at $\gamma(z)=\alpha_1,\ldots,\alpha_n$.

For a particular choice of one of the parameters and a discretization
of the values of $z$, the functions $R_n(z;q,p)$ yield
elliptic generalizations of Wilson's finite-dimensional
$_9F_8$ and $_{10}\Phi_9$ rational functions \cite{wil:orthogonal}.
They were derived in \cite{spi-zhe:spectral} from the theory of
self-similar solutions of nonlinear integrable discrete time chains
(for a brief review of the corresponding approach to special
functions, see \cite{spi:solitons,spi:factorization}).
Discrete analogs of equations (\ref{eig}),
(\ref{gevp}) valid for the latter finite-dimensional system of functions
were derived in \cite{spi-zhe:gevp} with the help of self-duality.
An equation satisfied by $_{10}\Phi_9$ functions, appearing
from $R_n(z;q,p)$ in the $p\to 0$ limit, was investigated
by Rahman and Suslov in \cite{rah-sus:classical}. The general three
term recurrence relations (\ref{ttr}) were considered in
\cite{zhe:bio} and, in a different form related to $R_{II}$ continued
fractions, in \cite{ism-mas:general}.

The solutions of the generalized eigenvalue problems are known to be
biorthogonal to each other; see, e.g.,
\cite{spi-zhe:spectral,spi-zhe:theory,zhe:bio} and the references therein.
Here we would like to demonstrate that the
elliptic beta integral (\ref{ell-int}) serves as the biorthogonality
measure for solutions of equation (\ref{gevp}). Consider
the scalar product
\begin{equation}
\int_C\Delta_E(z;\mathbf{t}) \Psi(z)\left(
\mathcal{D}_\eta-\lambda \mathcal{D}_\xi\right)\Phi(z)\frac{dz}{z},
\label{scalar}\end{equation}
where $\Delta_E(z;\mathbf{t})$ is the integrand of (\ref{ell-int})
and $\Phi(z), \Psi(z)$ are some complex functions.
The expression (\ref{scalar}) can be rewritten as
\begin{eqnarray}\nonumber
&& \int_{C}\Delta_E(z;\mathbf{t})(\kappa_\eta-V_\eta(z)-V_\eta(z^{-1}))
\Psi(z)\Phi(z)\frac{dz}{z}
\\ \nonumber && \makebox[2em]{}
+\int_{C_-}\Delta_E(q^{-1}z;\mathbf{t}) V_\eta(q^{-1}z)
\Psi(q^{-1}z)\Phi(z)\frac{dz}{z}
\\ && \makebox[2em]{}
+ \int_{C_+}\Delta_E(qz;\mathbf{t}) V_\eta(q^{-1}z^{-1})
\Psi(qz)\Phi(z)\frac{dz}{z} - \lambda \{\eta\to\xi\},
\label{scalar'}\end{eqnarray}
where $\{\eta\to\xi\}$ means the preceding expression with
$\eta$ replaced by $\xi$. The integration contours $C_\pm$ are obtained
from $C$ after the scaling transformations $z\to q^{\pm 1} z$. Suppose
that the poles of $\Delta_E(z;\mathbf{t})$ and the singularities of the
functions $\Phi(z), \Psi(q^{\pm 1}z)$ do not lie in the region swept by
the contours $C_\pm$ during their deformations to $C$.
Then (\ref{scalar'}) takes the form
\begin{equation}
\int_C\Delta_E(z;\mathbf{t}) \Phi(z)\left(
\mathcal{D}_\eta^T-\lambda \mathcal{D}_\xi^T\right)\Psi(z)\frac{dz}{z},
\label{scalar-conj}\end{equation}
where the adjoint (or transposed) operator $\mathcal{D}_\xi^T$
has the form
\begin{eqnarray} \nonumber
&& \mathcal{D}_\xi^T=\frac{\Delta_E(qz;\mathbf{t})}{\Delta_E(z;\mathbf{t})}
V_\xi(q^{-1}z^{-1})T+\frac{\Delta_E(q^{-1}z;\mathbf{t})}
{\Delta_E(z;\mathbf{t})}V_\xi(q^{-1}z)T^{-1}
\\ && \makebox[4em]{}
-V_\xi(z)-V_\xi(z^{-1})+\kappa_\xi.
\label{op-conj}\end{eqnarray}

Suppose that $\Phi_\lambda(z)$ is a solution of the equation
$\left(\mathcal{D}_\eta-\lambda \mathcal{D}_\xi\right)\Phi(z)=0$
and $\Psi_{\lambda'}(z)$ solves the conjugate equation
$\left(\mathcal{D}_\eta^T-\lambda'\mathcal{D}_\xi^T\right)\Psi(z)=0$
for some $\lambda'$. Both these functions can be multiplied by arbitrary
functions $f(z)$ satisfying the condition of periodicity on the logarithmic
scale, $f(qz)=f(z)$.  After the replacement of $\Phi(z)$ and $\Psi(z)$ in
(\ref{scalar}) and (\ref{scalar-conj}) by $\Phi_\lambda(z)$ and
$\Psi_{\lambda'}(z)$, these expressions become equal to zero.
In particular, (\ref{scalar-conj}) yields the relation
\begin{eqnarray}\nonumber
&& \int_C\Delta_E(z;\mathbf{t}) \Phi_\lambda(z)\left(
\mathcal{D}_\eta^T-\lambda \mathcal{D}_\xi^T\right)\Psi_{\lambda'}(z)
\frac{dz}{z}
\\ && \makebox[4em]{}
=(\lambda'-\lambda)\int_C\Delta_E(z;\mathbf{t})
\Phi_\lambda(z)\mathcal{D}_\xi^T\Psi_{\lambda'}(z)\frac{dz}{z}=0,
\label{biort}\end{eqnarray}
which shows that for $\lambda'\neq\lambda$ the function $\Phi_\lambda(z)$
is orthogonal to $\mathcal{D}_\xi^T\Psi_{\lambda'}(z).$

We find a function $g(z)$ such that
$$
g^{-1}(z)\left(\mathcal{D}_\eta^T-\lambda \mathcal{D}_\xi^T\right)g(z)
=\mathcal{D}_\eta-\lambda \mathcal{D}_\xi.
$$
After substitution of the known expressions for
$\Delta_E(qz;\mathbf{t})/\Delta_E(z;\mathbf{t})$ and the
spectral parameter (see $\lambda$ in (\ref{lambda})),
we get the following equation for $g(z)$:
$$
g(qz)=\frac{\theta\left(\frac{q}{t_4z},\frac{q\mu z}{t_4},
Az,\frac{A\mu}{q^2z};p\right)}{\theta\left(\frac{q^2z}{t_4},
\frac{\mu}{t_4 z},\frac{A}{qz},\frac{A\mu z}{q};p\right)}\, g(z),
$$
which is  solved easily:
\begin{equation}
g(z)=\frac{\Gamma(\frac{q\mu z}{t_4},\frac{\mu q}{t_4 z},
Az,\frac{A}{z};q,p)}{\Gamma(\frac{q^2z}{t_4},\frac{q^2}{t_4 z},
\frac{A\mu z}{q}, \frac{A\mu}{qz};q,p)}.
\label{g-function}\end{equation}
Here we have neglected the arbitrary factor $f(qz)=f(z)$, which has already
been mentioned. As a result, we get a direct relation between
$\Phi_\lambda(z)$ and $\Psi_\lambda(z)$:
$\Psi_\lambda(z)=g(z)\Phi_\lambda(z)$, where $g(z)$ depends on
$\lambda$ as well.

We denote $\lambda_n\equiv \lambda|_{\mu=q^n}$ and
$$
g_n(z)\equiv g(z)|_{\mu=q^n}=
\frac{\theta\left(\frac{q^2z}{t_4},\frac{q^2}{t_4z};p;q\right)_{n-1}}
{\theta\left(Az,\frac{A}{z};p;q\right)_{n-1}}.
$$
Substituting $\Phi_{\lambda_n}(z)=R_n(z;q,p)$ in the derived
biorthogonality relations, we see that the functions $R_n(z;q,p)$ are
formally orthogonal to $\mathcal{D}_\xi^T g_m(z)R_m(z;q,p)$ for $n\neq m$.

The general considerations of \cite{spi-zhe:theory,zhe:bio} show that
$R_n(z;q,p)$, which are rational functions of $\gamma(z)$ with the poles
at $\gamma(z)=\alpha_1,\ldots,\alpha_n$, are orthogonal to other rational
functions of $\gamma(z)$, which we denote as $T_m(z;q,p)$, with
the poles at $\gamma(z)=\beta_1,\ldots,\beta_{n}$. The choice of
$\alpha_n,\beta_n$ and the other recurrence coefficients in (\ref{ttr})
determine $R_n$ and $T_n$ uniquely, so that permutation of all $\alpha_n$
with $\beta_n$ permutes $R_n$ and $T_n$. In our case, we see
that parameters $\beta_n$ are obtained from
$\alpha_n=\gamma(q^n/t_4)$ after the replacement of $t_4$ by $pq/A$.
Equivalently, this replacement converts $\beta_n$ to $\alpha_n$.
An important point is that the weight function $\Delta_E(z,\mathbf{t})$
is invariant under such a transformation. Therefore, we can get $T_n$ out
of $R_n$ simply by the $t_4\to pq/A$  involution, which yields
\begin{equation}
T_n(z;q,p)={_{12}V_{11}}\left(\frac{At_3}{q};\frac{A}{t_0},\frac{A}{t_1},
\frac{A}{t_2},t_3z,\frac{t_3}{z},q^{-n},\frac{Aq^{n-1}}{t_4};q,p\right),
\label{T_n}\end{equation}
where the dependence on $p$ in the parameters drops out due to the
total ellipticity of this $_{12}V_{11}$ series.
Comparing with the previous consideration, we see that
$$
\mathcal{D}_\xi^T g_n(z)R_n(z;q,p) = \rho_n T_n(z;q,p)
$$
for some proportionality constants $\rho_n$, which are of no importance
 for us here.

Thus, the operator formalism developed above leads to the
following formal biorthogonality relation:
\begin{equation}
\int_C T_n(z;q,p)R_m(z;q,p)\Delta_E(z,\mathbf{t})\frac{dz}{z}=\tilde h_n
\delta_{nm}
\label{formal-ort}\end{equation}
for some constants $\tilde h_n$. Suppose that $C=\mathbb{T}$
and $|t_r|<1, |qp|<|A|$. Some poles of the functions $R_n,T_m$
cancel with zeros of $\Delta_E(z;\mathbf{t})$. The remaining ones
approach $\mathbb{T}$ as the indices $n,m$ increase. Starting with
sufficiently high values of $n$ and $m$, the contour $\mathbb{T}$ stops
to satisfy the conditions used in the derivation of (\ref{scalar-conj})
and must be deformed. As is shown in Appendix A,
(\ref{formal-ort}) is true if $C$ separates the points
$z= t_{0,1,2,3}p^aq^b,$ $t_4p^{a}q^{b-m},$ and
$A^{-1}p^{a+1}q^{b+1-n},$ $a,b\in\mathbb{N},$ from their partners
with the inverse $z\to z^{-1}$ coordinates.

Relation (\ref{formal-ort}) seems to remain true even if we multiply
$R_n(z;q,p)$ or $T_n(z;q,p)$ by an arbitrary function $f(z)$ with the
property $f(qz)=f(z)$. However, such nontrivial $f(z)$
must contain singularities which are crossed over when the contours
$C_\pm$ are deformed to $C$ (otherwise $f(z)=const$). Therefore, the
influence of such additional factors should be considered more
carefully. Moreover, only for very special $f(z)$
the normalization constants $\tilde h_n$ may admit exact evaluation.

The weight function $\Delta_E(z,\mathbf{t})$ is symmetric in $q$ and $p$,
whereas neither $R_n(z;q,p)$ nor $T_n(z;q,p)$ possess such a property.
We can try to restore this symmetry using the freedom in the factor
$f(z)=f(qz)$. Take $f(z)=R_k(z;p,q)$, $k\in\mathbb{N}$, that is, to the
functions $R_n(z;q,p)$ themselves with the permuted bases $q$ and $p$.
Then, the product
$$
R_{nk}(z)\equiv R_n(z;q,p)R_k(z;p,q)
$$
satisfies two generalized eigenvalue problems: (\ref{eig}) and the
$p$-difference equation obtained from it by the permutation of $q$
and $p$. For (\ref{eig}) we should have $\mu=q^n$, and for its
partner $\mu=p^k$. The function (\ref{lambda}) does not change under
the substitution $\mu\to p\mu$. Therefore, the choice
$\mu=q^np^k,\, n,k\in \mathbb{N},$
gives ``the spectrum" for both generalized eigenvalue problems.
The first factor of $R_{nk}(z)$ is a rational function of $\gamma(z;p)$
(we indicate the dependence on the base $p$ explicitly), but
the second is rational in $\gamma(z;q)$. Therefore, for generic $q,p$
it is necessary to view the functions $R_{nk}(z)$ not as rational functions
of some variable but as meromorphic functions of $z$.

In the same way,
the series termination condition $t_6=q^{-n}$ in (\ref{R_n}) may
be replaced by $t_6=q^{-n}p^{-k}$, which terminates simultaneously
the $_{12}V_{11}$ series for $R_k(z;p,q)$. The property of the total
ellipticity of the balanced $_{r+1}V_{r}(t_0;$ $t_1,\ldots,$ $t_{r-4};q,p)$
series plays a crucial role at this place: any parameter
$t_1,\ldots,t_{r-5}$ may be multiplied by an arbitrary integer power
of $p$ without changes (note that the parameter $t_0$ plays a
distinguished role and the series are invariant under the
transformation $t_0\to p^2t_0$).

As a result of the doubling of eigenvalue problems, the functions $R_{nk}(z)$
turn out to satisfy a quite unusual biorthogonality relation
(\ref{ort2}) characteristic
of functions of {\em two} independent variables, which was announced in
\cite{spi:special}. A rigorous consideration of all these biorthogonalities
associated with the function $R_n(z;q,p)$ with complete proofs is given
in Appendix A.

In fact, there is a deeper relationship between the structure of the
elliptic beta integral (\ref{ell-int}) and the biorthogonal functions
$R_n(z;q,p), T_n(z;q,p)$ and the solutions of equation (\ref{eig}) than
it is indicated in this section. We hope to address this later on.
As far as the multivariable generalizations
are concerned, there are some multidimensional analogs of equation
(\ref{eig}) for the solutions of which the multiple elliptic beta
integrals on root systems described in this paper may serve as
biorthogonality measures. However, their consideration lies beyond the
scope of the present work.

\appendix

\section{Proof of a biorthogonality relation}

We define a pair of functions of $z\in\mathbb{C}$ as products of two
terminating $_{12}V_{11}$ very-well-poised balanced theta hypergeometric
series with the argument $x=1$ and {\em different modular parameters}:
\begin{eqnarray} \nonumber
&&R_{nm}(z)={_{12}V_{11}}\left(\frac{t_3}{t_4};
\frac{q}{t_0t_4},\frac{q}{t_1t_4},\frac{q}{t_2t_4},t_3z,\frac{t_3}{z},q^{-n},
\frac{Aq^{n-1}}{t_4};q,p\right) \\  \label{R_nm}
&&\makebox[3em]{}\times {_{12}V_{11}}\left(\frac{t_3}{t_4};
\frac{p}{t_0t_4},\frac{p}{t_1t_4},\frac{p}{t_2t_4},t_3z,\frac{t_3}{z},p^{-m},
\frac{Ap^{m-1}}{t_4};p,q\right), \\ \nonumber
&&T_{nm}(z)={_{12}V_{11}}\left(\frac{At_3}{q};\frac{A}{t_0},\frac{A}{t_1},
\frac{A}{t_2},t_3z,\frac{t_3}{z},q^{-n},\frac{Aq^{n-1}}{t_4};q,p\right)
\\ \label{T_nm}
&& \makebox[3em]{}\times {_{12}V_{11}}\left(\frac{At_3}{p};\frac{A}{t_0},\frac{A}{t_1},
\frac{A}{t_2},t_3z,\frac{t_3}{z},p^{-m},\frac{Ap^{m-1}}{t_4};p,q\right),
\end{eqnarray}
where $n,m\in \mathbb{N}$. Obviously, these functions are symmetric
with respect to the permutation of $p$ and $q$. The balancing
conditions are used already in these series in order to
express one of the parameters in terms of the others.

For $m\neq 0$ and the fixed parameters $q, t_r, r=0,\ldots,4,$
the limit as $p\to 0$ is not well defined for $R_{nm}(z)$ and $T_{nm}(z)$.
The reason for this comes from the quasiperiodicity of
$\theta(z;p)$ because the limits $z\to 0$ or $z\to\infty$ are
not defined for it.

Now, we consider the following integral:
\begin{equation}\label{ort-int}
J_{mn,kl}\equiv\int_{C_{mn,kl}}T_{nl}(z)R_{mk}(z)
\Delta_E(z,\mathbf{t}) \frac{d z}{z},
\end{equation}
where $\Delta_E(z,\mathbf{t})$ denotes the weight function (\ref{weight})
and $C_{mn,kl}$ is some contour of integration.
We want to show that a particular choice of
$C_{mn,kl}$ (to be specified below) leads to the formula
\begin{equation}\label{ort}
J_{mn,kl}=h_{nl}\delta_{mn}\delta_{kl},
\end{equation}
where $h_{nl}$ are some normalization constants.

An elliptic extension of the Bailey transformation for the $_{10}\Phi_9$
series, which was derived by Frenkel and Turaev in
\cite{fre-tur:elliptic}, has the following form
(an alternative proof of it was given in \cite{spi:bailey}):
\begin{eqnarray}\nonumber
{_{12}V_{11}}(t_0;t_1,\dots,t_7;q,p) &=&
\frac{\theta(qt_0,qs_0/t_4,qs_0/t_5,qt_0/t_4t_5;p;q)_N}
{\theta(qs_0,qt_0/t_4,qt_0/t_5,qs_0/t_4t_5;p;q)_N}   \\
&& \times {_{12}V_{11}}(s_0;s_1,\dots,s_7;q,p),
\label{ft-bailey}\end{eqnarray}
where $\prod_{m=1}^7t_m=t_0^3q^2$, $t_6=q^{-N}$, $N\in\mathbb{N}$,
\begin{equation}
s_0=\frac{qt_0^2}{t_1t_2t_3},\quad s_1=\frac{s_0t_1}{t_0},\quad
s_2=\frac{s_0t_2}{t_0},\quad s_3=\frac{s_0t_3}{t_0},
\label{s-param}\end{equation}
and $s_4, s_5, s_6, s_7$ form an arbitrary permutation of
$t_4, t_5, t_6, t_7$. Using this identity,
we can rewrite the functions $R_{mk}(z)$ as follows
\begin{eqnarray} \nonumber
&& R_{mk}(z) = \kappa_m(p;q)\kappa_k(q;p) \\
&& \makebox[3em]{}  \times
{_{12}V_{11}}\left(\frac{t_0}{t_4};\frac{q}{t_1t_4},\frac{q}{t_2t_4},
\frac{q}{t_3t_4},t_0z,\frac{t_0}{z},q^{-m},\frac{Aq^{m-1}}{t_4};q,p\right)
\nonumber \\
&& \makebox[3em]{} \times
{_{12}V_{11}}\left(\frac{t_0}{t_4};\frac{p}{t_1t_4},\frac{p}{t_2t_4},
\frac{p}{t_3t_4},t_0z,\frac{t_0}{z},p^{-k},\frac{Ap^{k-1}}{t_4};p,q\right),
\label{R-transf}\end{eqnarray}
where
$$ \kappa_m(p;q)=\frac{\theta(qt_3/t_4,t_0t_1,t_0t_2,A/qt_0;p;q)_m}
{\theta(qt_0/t_4,t_1t_3,t_2t_3,A/qt_3;p;q)_m}. $$

Substituting the explicit series expressions \eqref{T_nm} and
\eqref{R-transf} in the definition of $J_{mn,kl}$ (\ref{ort-int}) yields
\begin{eqnarray}\nonumber
J_{mn,kl}&=& \kappa_m(p;q)\kappa_{k}(q;p)
\sum_{r=0}^n\sum_{r'=0}^l\sum_{s=0}^{m}\sum_{s'=0}^k q^{r+s}p^{r'+s'}
\\ \nonumber
&& \times \frac{\theta(At_3q^{2r-1},t_0q^{2s}/t_4;p)}{\theta(At_3/q,t_0/t_4;p)}
\frac{\theta(At_3p^{2r'-1},t_0p^{2s'}/t_4;q)}{\theta(At_3/p,t_0/t_4;q)}
\\  \nonumber
&& \times \frac{\theta(At_3/q,A/t_0,A/t_1,A/t_2,q^{-n},Aq^{n-1}/t_4;p;q)_r}
{\theta(q,t_0t_3,t_1t_3,t_2t_3,At_3q^n,t_3t_4q^{1-n};p;q)_r}
\\ \nonumber
&& \times \frac{\theta(At_3/p,A/t_0,A/t_1,A/t_2,p^{-l},Ap^{l-1}/t_4;q;p)_{r'}}
{\theta(p,t_0t_3,t_1t_3,t_2t_3,At_3p^l,t_3t_4p^{1-l};q;p)_{r'}}
\\ \nonumber
&& \times\frac{\theta(t_0/t_4,q/t_1t_4,q/t_2t_4,q/t_3t_4,q^{-m},
Aq^{m-1}/t_4;p;q)_s}
{\theta(q,t_0t_1,t_0t_2,t_0t_3,t_0q^{m+1}/t_4,t_0q^{2-m}/A;p;q)_s}
\\ \nonumber
&& \times\frac{\theta(t_0/t_4,p/t_1t_4,p/t_2t_4,p/t_3t_4,p^{-k},
Ap^{k-1}/t_4;q;p)_{s'}}
{\theta(p,t_0t_1,t_0t_2,t_0t_3,t_0p^{k+1}/t_4,t_0p^{2-k}/A;q;p)_{s'}}\,
I_{rs,r's'},
\end{eqnarray}
where
\begin{eqnarray}\nonumber
&I_{rs,r's'}& = \int_{C_{mn,kl}}\Delta_E(z,\mathbf{t})
\frac{\theta(zt_3,z^{-1}t_3;p;q)_r}{\theta(zA,z^{-1}A;p;q)_r}
\frac{\theta(zt_0,z^{-1}t_0;p;q)_s}
{\theta(zqt_4^{-1},z^{-1}qt_4^{-1};p;q)_s} \\
&& \times
\frac{\theta(zt_3,z^{-1}t_3;q;p)_{r'}}{\theta(zA,z^{-1}A;q;p)_{r'}}
\frac{\theta(zt_0,z^{-1}t_0;q;p)_{s'}}
{\theta(zpt_4^{-1},z^{-1}pt_4^{-1};q;p)_{s'}} \frac{d z}{z}.
\label{int3}\end{eqnarray}
We introduce the notation
$$
\tilde t_0=t_0q^sp^{s'},\quad \tilde t_1=t_1,\quad  \tilde t_2=t_2,
\quad \tilde t_3=t_3q^rp^{r'}, \quad \tilde t_4=t_4q^{-s}p^{-s'},
$$
so that
$$
\tilde A=\prod_{m=0}^4\tilde t_m=Aq^rp^{r'}.
$$
Then, using the transformation property
$$
\theta(z;p;q)_l=(-z)^lq^{l(l-1)/2}\theta(z^{-1}q^{-l+1};p;q)_l
=\frac{(-z)^lq^{l(l-1)/2}}{\theta(qz^{-1};p;q)_{-l}},
$$
we can rewrite the integral (\ref{int3}) in the form
\begin{equation}\label{int4}
I_{rs,r's'}=\left(\frac{t_0t_4}{pq}\right)^{2ss'}
\left(\frac{t_3}{A}\right)^{2rr'}\frac{t_4^{2(s+s')}}{q^{s(s+1)}p^{s'(s'+1)}}
\int_{C_{mn,kl}}
\Delta_E(z,\mathbf{\tilde t}) \frac{d z}{z}.
\end{equation}
But the integral on the right-hand side of (\ref{int4}) coincides
with the elliptic beta integral (\ref{ell-int}), provided we identify
$C_{mn,kl}=\mathbb{T}$ and impose the constraints $|\tilde t_m|<1$,
$|pq|<|\tilde A|$. However, the values of the integers
$r,s,r',s'\in\mathbb{N}$ are not limited; starting with their
sufficiently large values, we shall have either
$|\tilde A|=|q^rp^{r'}A|<|pq|$ or $|\tilde t_4|=|q^{-s}p^{-s'}t_4|>1$.
Now is the moment to specify the contour $C_{mn,kl}$. We choose
it in such a way that formula (\ref{ell-int}) remains applicable.
More precisely, let $C_{mn,kl}$ be a deformation of $\mathbb{T}$ such that
it separates the poles at $z=t_{0,1,2,3}p^aq^{b},$ $t_4p^{a-k}q^{b-m},$
and $A^{-1}p^{a+1-l}q^{b+1-n},$ $a,b\in\mathbb{N}$, that lie inside
$C_{mn,kl}$ and converge to zero,
from the poles diverging to infinity, which are obtained from the
poles inside $C_{mn,kl}$ by the inversion transformation $z\to z^{-1}$.
The subscripts $m,n,k,l$ in
the notation for such a contour indicate that, evidently, the
shape of $C_{nm,kl}$ depends on the indices of the functions
$T_{nl}(z), R_{mk}(z)$.

For such a contour $C_{mn,kl}$, we have
\begin{eqnarray}\nonumber
I_{rs,r's'} &=&
\left(\frac{t_0t_4}{pq}\right)^{2ss'}
\left(\frac{t_3}{A}\right)^{2rr'}\frac{t_4^{2(s+s')}}{q^{s(s+1)}p^{s'(s'+1)}}
\mathcal{N}_E(\mathbf{\tilde t}) \\  \nonumber
&=&\mathcal{N}_E(\mathbf{t})
\frac{\theta(t_1t_3,t_2t_3,t_3t_4;p;q)_r}
{\theta(A/t_0,A/t_1,A/t_2;p;q)_r}  \\  \nonumber
&&\times \frac{\theta(t_0t_3;p;q)_{r+s}}{\theta(A/t_4;p;q)_{r+s}}
\frac{\theta(t_0t_1,t_0t_2,q^{1-r}t_0/A;p;q)_s}
{\theta(q/t_1t_4,q/t_2t_4,q^{1-r}/t_3t_4;p;q)_s}
\\  \nonumber
&& \times
\frac{\theta(t_1t_3,t_2t_3,t_3t_4;q;p)_{r'}}
{\theta(A/t_0,A/t_1,A/t_2;q;p)_{r'}}  \\  \nonumber
&&\times \frac{\theta(t_0t_3;q;p)_{r'+s'}}{\theta(A/t_4;q;p)_{r'+s'}}
\frac{\theta(t_0t_1,t_0t_2,p^{1-r'}t_0/A;q;p)_{s'}}
{\theta(p/t_1t_4,p/t_2t_4,p^{1-r'}/t_3t_4;q;p)_{s'}},
\end{eqnarray}
where the function $\mathcal{N}_E(\mathbf{t})$ is fixed in \eqref{result}.

As a result of these manipulations, the quantity $J_{mn,kl}$
splits into a product of two double series each depending only on
the indices $m, n$ and $k,l$ separately. After an application of
the relation $(a;p;q)_{r+s}=(aq^r;p;q)_s(a;p;q)_r$ and various
simplifications, we can write
$$
J_{mn,kl}=\mathcal{N}_E(\mathbf{t})J_{mn}(p;q)J_{kl}(q;p),
$$
where
\begin{eqnarray}\nonumber
J_{mn}(p;q) &=& \kappa_m(p;q) \sum_{r=0}^nq^r
\frac{\theta(At_3q^{2r-1};p)}{\theta(At_3/q;p)}
\frac{\theta(At_3/q,q^{-n},Aq^{n-1}/t_4,t_3t_4;p;q)_r}
{\theta(q,At_3q^n,t_3t_4q^{1-n},A/t_4;p;q)_r}
\\  \label{8E7-aux}
&& \times {_{10}V_9}\left(\frac{t_0}{t_4};\frac{q}{t_3t_4},
t_0t_3q^r,\frac{q^{1-r}t_0}{A},\frac{Aq^{m-1}}{t_4},q^{-m};q,p\right).
\end{eqnarray}

The constraint $t_2t_3=qt_0$ imposed upon relation (\ref{ft-bailey})
converts the $_{12}V_{11}$-series on its left-hand side into a
terminating $_{10}V_9$
series, whereas on the right-hand side only the first term of the
corresponding $_{12}V_{11}$-series survives. As a result, we get the
Frenkel-Turaev sum, or an elliptic generalization of the Jackson sum
for terminating very-well-poised balanced $_8\Phi_7$-series:
\begin{eqnarray}\nonumber
\lefteqn{ {_{10}V_9}(t_0;t_1,t_4, t_5,t_6,t_7;q,p) } &&
\\ && \makebox[4em]{}
= \frac{\theta (qt_0,qt_0/t_1t_4,qt_0/t_1t_5,qt_0/t_4t_5;p;q)_N}
{\theta(qt_0/t_1t_4t_5,qt_0/t_1,qt_0/t_4,qt_0/t_5;p;q)_N},
\label{ft-sum}\end{eqnarray}
where $t_1t_4t_5t_6t_7 =qt_0^2$ and $t_6=q^{-N}$, $N\in \mathbb{N}$.
Application of this sum to the $_{10}V_9$ series in \eqref{8E7-aux}
yields
$$
{_{10}V_9}(\ldots)=\frac{\theta(qt_0/t_4,t_1t_2,At_3q^{r-1},q^{-r};p;q)_m}
{\theta(t_0t_3,A/qt_0,Aq^r/t_4,q^{1-r}/t_3t_4;p;q)_m}.
$$
Clearly, this expression vanishes for $m>r$. This means that $J_{mn}=0$
for $m>n$. For $m\leq n$, we get
\begin{eqnarray}\nonumber
J_{mn}(p;q)=\kappa_m(p;q)
\frac{\theta(At_3;p;q)_{2m}}{\theta(A/t_4;p;q)_{2m}}
\frac{\theta(Aq^{n-1}/t_4,qt_0/t_4,t_1t_2,q^{-n};p;q)_m}
{\theta(t_3t_4q^{1-n},t_0t_3,A/qt_0,At_3q^n;p;q)_m}
\\  \nonumber
\times (t_3t_4)^m\, {_8V_7}(At_3q^{2m-1};t_3t_4,Aq^{n+m-1}/t_4,q^{m-n};q,p).
\end{eqnarray}
Applying the summation formula (\ref{ft-sum}) to the latter
$_8V_7$ series, we get
$$
{_8V_7}(\ldots)=\frac{\theta(At_3q^{2m},q^{m-n+1},Aq^{m+n}/t_4,qt_3t_4;p;q)_{n-m}}
{\theta(q,Aq^{2m}/t_4,t_3t_4q^{m-n+1},At_3q^{m+n};p;q)_{n-m}},
$$
which is equal to zero for $n>m$ due to the factor
$\theta(q^{m-n+1};p;q)_{n-m}$. As a result,
$J_{mn}(p;q)=h_n(p;q)\delta_{mn}$, where the normalization constants
have the form
\begin{equation}\label{norm}
h_n(p;q)=\frac{\theta(A/qt_4;p)
\theta(q,qt_3/t_4,t_0t_1,t_0t_2,t_1t_2,At_3;p;q)_nq^{-n}}
{\theta(Aq^{2n-1}/t_4;p)
\theta(1/t_3t_4,t_0t_3,t_1t_3,t_2t_3,A/qt_3,A/qt_4;p;q)_n}.
\end{equation}

The fact that $J_{mn}=0$ for $n\neq m$ provides
the desired biorthogonality relation (\ref{ort}).
We summarize the result obtained in the form of the theorem
that was announced in \cite{spi:special} (it is necessary to
apply the series notation introduced in \cite{spi:theta,spi:bailey}
and permute the parameters $t_3$ and $t_4$ in \cite{spi:special} in
order to match with the current presentation).

\begin{theorem}
Let $t_m,$ $\Delta_E(z,{\bf t}),$ $\mathcal{N}_E({\bf t})$ be
the same as in Theorem 1. Let $C_{mn,kl}$ denote
a positively oriented contour separating the points
$z=\{ t_{0,1,2,3}p^aq^b,$ $t_4p^{a-k}q^{b-m},$
$A^{-1}p^{a+1-l}q^{b+1-n}\}_{a,b\in\mathbb{N}}$
from the points with the inverse ($z\to z^{-1}$) coordinates.
Then $R_{mk}(z)$ and $T_{nl}(z)$ satisfy the
following biorthogonality relation
\begin{equation}\label{ort2}
\int_{C_{mn,kl}}T_{nl}(z)R_{mk}(z)
\Delta_E(z,{\bf t}) \frac{d z}{z}=
h_{nl}\mathcal{N}_E({\bf t})\delta_{mn}\delta_{kl},
\end{equation}
where $h_{nl}$ are the normalization constants,
\begin{eqnarray} \nonumber
&& h_{nl}= \frac{\theta(A/qt_4;p)
\theta(q,qt_3/t_4,t_0t_1,t_0t_2,t_1t_2,At_3;p;q)_nq^{-n}}
{\theta(Aq^{2n}/qt_4;p)
\theta(1/t_3t_4,t_0t_3,t_1t_3,t_2t_3,A/qt_3,A/qt_4;p;q)_n} \\
&&\makebox[2em]{}\times  \frac{\theta(A/pt_4;q)
\theta(p,pt_3/t_4,t_0t_1,t_0t_2,t_1t_2,At_3;q;p)_lp^{-l}}
{\theta(Ap^{2l}/pt_4;q)
\theta(1/t_3t_4,t_0t_3,t_1t_3,t_2t_3,A/pt_3,A/pt_4;q;p)_l}.
\label{norm2}\end{eqnarray}
\end{theorem}

As is clear from (\ref{R_n}) and (\ref{T_n}), we have
$R_m(z;q,p)$ $= R_{m0}(z)$ and $T_n(z;q,p)= T_{n0}(z)$.
These functions $R_m, T_n$ are equal to $_{12}V_{11}$ elliptic
hypergeometric series with particular choices of the parameters.

\begin{corollary}
The functions $R_m(z;q,p)$ and $T_n(z;q,p)$ satisfy the following
bi\-or\-tho\-gonality condition
\begin{equation}\label{R_mT_n-ort}
\int_{C_{mn}}T_{n}(z;q,p)R_{m}(z;q,p)
\Delta_E(z,{\bf t}) \frac{d z}{z}= h_{n}\mathcal{N}_E({\bf t})\delta_{mn},
\end{equation}
where the constants $h_n$ are fixed in (\ref{norm})
and the contour $C_{mn}$ encircles the poles of the integrand located
at $z=\{ t_{0,1,2,3}q^ap^b,
t_4p^aq^{b-m},A^{-1}p^{a+1}q^{b+1-n}\}_{a,b\in\mathbb{N}}$ and separates
them from the poles with inverse $z\to z^{-1}$ coordinates.
\end{corollary}

The biorthogonal rational functions $R_m(z;q,p)$ and $T_n(z;q,p)$
describe elliptic generalizations of the Rahman set of
continuous $_{10}\Phi_9$ functions \cite{rah:integral}
to which they are reduced in the limit as $p\to 0$. Accordingly,
in this limit, formula (\ref{R_mT_n-ort}) is reduced to the Rahman
biorthogonality condition.

\section{Integral representations for $_{12}E_{11}$ series}

Here we derive an integral representation for the product of two
terminating $_{12}E_{11}$ (more precisely, $_{12}V_{11}$) series with
some particular choice of parameters. For this, we apply an elliptic
generalization of the technique
used in \cite{rah:integral} for the derivation of the contour integral
representation for a terminating $_{10}\Phi_9$ series.

\begin{theorem}
Suppose that five parameters $t_k, k=0,\ldots, 4,$ satisfy the conditions
of Theorem 1. Denote by $m,n$ two positive integers and by $C_{mn}$
a positively oriented contour such that for all $a,b \in\mathbb{N}$
it separates the points $z=\{t_kp^aq^b,$ $A^{-1}q^{b+1-m}p^{a+1-n}\}$
from their partners with the inverse coordinates ($z\to z^{-1}$).
Under these conditions, the following integral representation for the
product of two $_{12}V_{11}$ terminating very-well-poised
balanced theta hypergeometric series at $x=1$ holds true:
\begin{eqnarray}\nonumber
&&{_{12}V_{11}}\left(\frac{At_0}{q};\alpha,t_0t_1,t_0t_2,t_0t_3,t_0t_4,
q^{-m},\frac{A^2q^{m-1}}{\alpha};q,p\right) \\  \nonumber
&& \times {_{12}V_{11}}
\left(\frac{At_0}{p};\beta,t_0t_1,t_0t_2,t_0t_3,t_0t_4,
p^{-n},\frac{A^2p^{n-1}}{\beta};p,q\right) \\  \nonumber
&&=\frac{1}{\mathcal{N}_E({\bf t})}
\frac{\theta(At_0,\frac{A}{t_0};p;q)_m\theta(At_0,\frac{A}{t_0};q;p)_n}
{\theta(\frac{A}{\alpha t_0},\frac{At_0}{\alpha};p;q)_m
\theta(\frac{A}{\beta t_0},\frac{At_0}{\beta};q;p)_n}  \\
&& \makebox[1em]{}
\times \int_{C_{mn}}
\Delta_E(z,{\bf t})\frac{\theta(\frac{Az}{\alpha},\frac{A}{\alpha z};p;q)_m
\theta(\frac{Az}{\beta},\frac{A}{\beta z};q;p)_n }
{\theta(Az,\frac{A}{z};p;q)_m\theta(Az,\frac{A}{z};q;p)_n}\frac{d z}{z},
\label{intrep}\end{eqnarray}
where $\alpha$ and $\beta$ are arbitrary complex parameters.
\end{theorem}

\begin{proof}
Under the conditions imposed upon the parameters
in the formulation of this theorem, the following relations are true:
\begin{eqnarray}\nonumber
\lefteqn{
\int_{C_{ij}}\Delta_E(z,\mathbf{t})
\frac{\theta(zt_0,z^{-1}t_0;p;q)_i\theta(zt_0,z^{-1}t_0;q;p)_j}
{\theta(zA,z^{-1}A;p;q)_i\theta(zA,z^{-1}A;q;p)_j}\frac{dz}{z}
} && \\ &&
= \left(\frac{t_0}{A}\right)^{2ij}
\mathcal{N}_E(t_0q^ip^j,t_1,\ldots,t_4)
\nonumber \\  &&
=\frac{\theta(t_0t_1,\ldots,t_0t_4;p;q)_i\theta(t_0t_1,\ldots,t_0t_4;q;p)_j}
{\theta(A/t_1,\ldots,A/t_4;p;q)_i\theta(A/t_1,\ldots,A/t_4;q;p)_j}
\, \mathcal{N}_E(\mathbf{t}).
\label{identity}\end{eqnarray}

We multiply (\ref{identity}) by the factor
\begin{eqnarray*}
  && q^i\frac{\theta(At_0q^{2i-1};p)}{\theta(At_0/q;p)}
\frac{\theta(At_0/q,\alpha,q^{-m},A^2q^{m-1}/\alpha;p;q)_i}
     {\theta(q,At_0/\alpha,At_0q^m,\alpha q^{1-m}t_0/A;p;q)_i}
\\
 && \times p^j\frac{\theta(At_0p^{2j-1};q)}{\theta(At_0/p;q)}
\frac{\theta(At_0/p,\beta,p^{-n},A^2p^{n-1}/\beta;q;p)_j}
     {\theta(p,At_0/\beta,At_0p^n,\beta p^{1-n}t_0/A;q;p)_j},
\end{eqnarray*}
where $\alpha$ and $\beta$ are arbitrary complex parameters,
and sum over $i$ from 0 to $m$ and over $j$ from 0 to $n$. As a result,
we get the relation
\begin{eqnarray*}
&&\int_{C_{mn}} \Delta_E(z,\mathbf{t})\,
{_{10}V_9}(At_0/q;t_0z,t_0z^{-1},\alpha,q^{-m},A^2q^{m-1}/\alpha;q,p)
\\ && \makebox[4em]{}
\times {_{10}V_9}(At_0/p;t_0z,t_0z^{-1},\beta,p^{-n},A^2p^{n-1}/\alpha;p,q)
\,\frac{dz}{z}  \\
&&\makebox[2em]{} ={_{12}V_{11}}(At_0/q;\alpha,q^{-m},A^2q^{m-1}/\alpha,
t_0t_1,\ldots,t_0t_4;q,p) \\
&& \makebox[3em]{} \times
{_{12}V_{11}}(At_0/p;\beta,p^{-n},A^2p^{n-1}/\beta,t_0t_1,\ldots,t_0t_4;
p,q)\, \mathcal{N}_E(\mathbf{t}).
\end{eqnarray*}
Application of the Frenkel-Turaev sum to the $_{10}V_9$ series standing
under the integral sign leads to (\ref{intrep}).
\end{proof}

For $n=0$ we get an integral representation of a single terminating
$_{12}V_{11}$ series, which can be reduced further to the $_{10}\Phi_9$
$q$-series level \cite{rah:integral} by letting $p\to 0$. However, for
$n\neq 0$ the limit as $p\to 0$ is not well defined, and formula
(\ref{intrep}) exists only at the elliptic level.

\bibliographystyle{amsplain}

\end{document}